\RequirePackage{etoolbox}
\csdef{input@path}{{style/}{graphics/}}
\documentclass[aap,MSNbibl,nameyear,dvips]{arximspdf}
\makeatletter
   \@ifpackageloaded{graphicx}{}{\usepackage{graphicx}}
\makeatother

\doi{10.1214/14-AAP1017} 
\volume{25}
\issue{3}
\pubyear{2015}
\firstpage{1079}
\lastpage{1107}

\makeatletter
\newcommand{\rrvert}{\vert}
\newcommand{\llvert}{\vert}
\newcommand{\eqref}[1]{(\ref{#1})}
\newproclaim{definition}{Definition}[section]
\newtheorem{lemme}[definition]{Lemma}
\newtheorem{theorem}[definition]{Theorem}

\newproclaim{rem}[definition]{Remark}
\newtheorem{prop}[definition]{Proposition}
\newproclaim{propriete}[definition]{Property}

\def\gap{\operatorname{gap}}
\def\L{\mathcal{L}}
\def\Z{\mathbb{Z}}
\def\N{\mathbb{N}}
\def\R{\mathbb{R}}
\makeatother

\begin{document}
\begin{frontmatter}

\title{Tracer diffusion at low temperature in kinetically
constrained models}
\runtitle{Diffusion in KCSM}

\begin{aug}
\author[A]{\fnms{Oriane}~\snm{Blondel}\corref{}\ead[label=e1]{oriane.blondel@ens.fr}\ead[label=u1,url]{http://www.proba.jussieu.fr/pageperso/blondel/}}
\affiliation{Universit\'e Paris Diderot}
\address[A]{LPMA\\
Universit\'e Paris Diderot\\
Case courrier 7012\\
Avenue de France\\
75205 Paris Cedex 13\\
France\\
\printead{e1}\\
\printead{u1}}
\runauthor{O. Blondel}
\end{aug}

\received{\smonth{6} \syear{2013}}

%
\begin{abstract}
We describe the motion of a tracer in an environment given by a
kinetically constrained spin model (KCSM) at equilibrium. We check
convergence of its trajectory properly rescaled to a Brownian motion
and positivity of the diffusion coefficient $D$ as soon as the spectral
gap of the environment is positive (which coincides with the ergodicity
region under general conditions). Then we study the asymptotic behavior
of $D$ when the density $1-q$ of the environment goes to $1$ in two
classes of KCSM. For noncooperative models, the diffusion coefficient
$D$ scales like a power of $q$, with an exponent that we compute
explicitly. In the case of the Fredrickson--Andersen one-spin
facilitated model, this proves a prediction made in Jung, Garrahan and
Chandler [\textit{Phys. Rev. E} \textbf{69} (2004) 061205]. For the
East model,
instead we prove that the diffusion coefficient is comparable to the
spectral gap, which goes to zero faster than any power of $q$. This
result contradicts the prediction of physicists (Jung, Garrahan and
Chandler [\textit{Phys. Rev. E} \textbf{69} (2004) 061205;
\textit{J. Chem. Phys.} \textbf{123} (2005) 084509]),
based on numerical simulations, that
suggested $D\sim\gap^\xi$ with $\xi<1$.
\end{abstract}

%


\end{frontmatter}

\section{Introduction}
Kinetically constrained models (KCSM) have been introduced in the
physics literature to model glassy dynamics. They are Markov processes
on $\lbrace0,1\rbrace^{\Z^d}$ (or more generally on the set of
configurations on a graph), where zeros mark empty sites, and ones mark
sites occupied by a particle. The dynamics is of Glauber type: with
rate one, each site refreshes its occupation variable: to a zero with
probability $q$, and to a one with probability $1-q$, on the condition
that a specific constraint be satisfied by the configuration around the
to-be-updated site. This constraint takes the form that ``a certain set
of zeros should be present in a fixed neighborhood,'' but does not
involve the configuration \emph{at} the to-be-updated site, so that the
product Bernoulli measure on $\Z^d$ with parameter $1-q$ is reversible
for the dynamics.

A tracer particle evolves in an environment given by a KCSM. The
environment is not influenced by the tracer, which performs a simple
random walk constrained to jumping only between two empty sites.
Properly rescaled, the tracer trajectory is expected to converge to a
Brownian motion with a diffusion coefficient depending on the
environment. Standard results and strategy [\citet{kipnisvaradhan},
\citet
{demasiferrarigoldsteinwick}, \citet{spohn}] allow us to show that in
the ergodic regime for the environment there is indeed convergence to a
Brownian motion, and to give a variational formula for the diffusion
coefficient; see Proposition~\ref{prop:diffusion} and Lemma~\ref
{lem:varformula}. A general argument then implies that, as soon as the
environment has a positive spectral gap, the diffusion coefficient is
also positive, so that the convergence result is nondegenerate
(Proposition~\ref{prop:easybounds}). Note that the ergodicity regime of
KCSM has been identified in \citet{KCSM}, and has been shown to coincide
with the region of positivity of the spectral gap in great generality,
including all the models we consider. Thus we prove in fact positivity
of the diffusion coefficient in the ergodic regime of the dynamical
environment. The variational formula also yields an immediate upper
bound on the diffusion coefficient. A similar study was carried in
\citet
{bertinitoninelli} with environments given by some noncooperative
constrained models with Kawasaki dynamics.

The main focus of this paper is to compute the asymptotics of the
diffusion coefficient when $q\rightarrow0$. This study is inspired by
the papers Jung, Garrahan and
Chandler (\citeyear{junggarrahanchandler,junggarrahanchandler2}), which in turn have the following physical
motivation. In homogeneous liquid systems, physicists argue that the
relaxation time $\tau$ (measured as the viscosity of the liquid), the
temperature $T$ and the diffusion coefficient $D$ of a particle moving
inside the system satisfy the following relation, called the
Stokes--Einstein relation,
%
%
\begin{equation}
D\propto T\tau^{-1}.
\end{equation}

This relation is well obeyed in liquids at high enough temperature.
Instead, in supercooled liquids it is experimentally observed [see, for
instance, \citet{EEHPW}, \citet{CiceroneEdiger96}, \citet
{ChangSillescu97}, \citet{Swallenetal03}] that $D\tau/T$
increases by 2--3 orders of magnitude when decreasing $T$ toward the
glass transition temperature. In particular both $D$ and $\tau^ {-1}$
decrease faster than any power law when the temperature is lowered, and for
many supercooled liquids a good fit of data is
%
%
\begin{equation}
D\propto\tau^{-\xi} \qquad\mbox{with }\xi<1.
\end{equation}
In other words, the self-diffusion of particles becomes much faster than
structural relaxation, and the Stokes--Einstein relation is violated. This
decoupling between translational diffusion and global relaxation is
interpreted as a landmark of dynamical heterogeneities in glassy systems,
namely the existence of spatially correlated regions of relatively high or
low mobility that persist for a finite lifetime in the liquid, and that
grow in size as one approaches the glass transition. More precisely, the
decoupling should be due to the fact that diffusion is dominated by the
fastest regions, whereas structural relaxation is dominated by the slowest
regions.

In order to investigate the possible violation of the Stokes--Einstein
relation in KCM, which are used as simplified models of glassy
dynamics, in Jung, Garrahan and
Chandler (\citeyear{junggarrahanchandler,junggarrahanchandler2}) the authors run simulations of a tracer in two systems
with constrained dynamics in one dimension: the FA-1f model (in which
the constraint requests that at least one neighbor be empty) and the
East model (in which the constraint is satisfied if the neighbor in the
East direction is empty). They predict in both cases a breakdown of the
Stokes--Einstein relation. More precisely, they predict that in the
FA-1f model in one dimension,
%
%
\begin{equation}
\label{eq:predictionFA} D\sim q^2\sim\gap^{2/3}
\end{equation}
and in the East model,
%
%
\begin{equation}
\label{eq:predictioneast} D\approx\gap^{\xi} \qquad\mbox{with } \xi \approx0.73.
\end{equation}

Our results confirm \eqref{eq:predictionFA} but invalidate \eqref
{eq:predictioneast}. Indeed we prove that for the
East model $D \approx\gap$ up to polynomial corrections (Theorem~\ref
{th:kzeros}). For this
model simulations are much harder to run than for FA model due to the very
fast divergence of the relaxation time when $q \rightarrow0$ [faster
than any power of
$1/q$; see \eqref{eq:gapeast}], thus accounting for the wrong numerical
prediction.

More generally we show that, in any dimension, if the model is defined
by the constraint ``there should be at least $k$ zeros in a ball of
radius $k$ around the to-be-updated site,'' the diffusion coefficient is
of order $q^{k+1}$ [$k=1$ corresponds to the FA-1f model, so the result
confirms the conjecture in \citet{junggarrahanchandler}; see
Theorem~\ref
{th:kzeros}]. The proof of this result relies on the introduction of an
auxiliary dynamics whose diffusion coefficient gives a lower bound for
$D$. This dynamics is similar to that in \citet{spohn}, though it is
less immediate to derive because it does not appear by just suppressing
terms in the variational formula. The very construction of this
auxiliary dynamics is in fact quite informative about the effective
dynamics of the tracer, and can be generalized to other noncooperative
models; see Definition~\ref{def:noncoop}. Back to the FA-1f model, in
dimension 2, our result and the estimate of the spectral gap in \citet
{KCSM} (Theorem~6.4) show that $D\propto\gap$. When $d\geq3$, our
bounds allow us to extract the asymptotic dependence of $D$ in $q$.
However, due to the current lack of precise bounds on the spectral gap,
we cannot decide whether $D\propto\gap^\xi$ for some exponent $\xi$,
but our results do imply that $\xi$ cannot be strictly smaller than one.

We also study the diffusion coefficient when the environment is given
by the East model, which does not belong to the noncooperative class.
As mentioned above, we prove in this case $D \approx\gap$ up to polynomial
corrections (Theorem~\ref{th:east}), contradicting \eqref
{eq:predictioneast}. The strategy used in that context is very
different from the one we designed for the ``$k$-zeros'' model because
the dynamics of the East model is cooperative, so that restricting the
dynamics only to a neighborhood of the tracer is not relevant. The
proof relies instead on precise estimates of the energy barriers that
have to be overcome in order for the tracer to cross the typical
distance between two zeros at equilibrium, $1/q$. These estimates have
been established mostly in \citet{KCSM} and \citet{timescaleseparation}.
As an extension of results in these two papers, we provide in
particular a better estimate on the spectral gap in infinite volume
(Lemma~\ref{lem:newgapeast}).

The paper is organized as follows. In Section~\ref{sec:models}, we
define the processes of the environment, the tracer dynamics and the
environment seen from the tracer. In Section~\ref{sec:results} we
collect the main results of this paper, which are proved in the
following sections. In Section~\ref{sec:diffusivity}, we prove
convergence of the tracer trajectory to a Brownian motion with positive
diffusion coefficient in the ergodic regime. Section~\ref{sec:kzeros}
is devoted to retrieve the right asymptotics for the diffusion
coefficient when the density goes to $1$ in noncooperative models.
Finally, in Section~\ref{sec:east}, we show that asymptotically the
diffusion coefficient in the East model is of the same order as the
spectral gap, up to polynomial corrections.

\section{Models and notation}\label{sec:models}
Let $\Omega=\lbrace0,1\rbrace^{\Z^d}$. For $\omega\in\Omega$,
$x\in\Z
^d$ we define $\omega^x$ the configuration such that
%
%
\begin{equation}
\omega^x_y = \cases{ %
\omega_y,& \quad $\mbox{if $y\neq x$},$
\vspace*{2pt}\cr
1-\omega_x, & \quad $\mbox{if $y=x$}$.}
\end{equation}
A KCSM is defined by its equilibrium density $p=1-q$ and constraints
$ (c_x(\omega) )_{x\in\Z,\omega\in\Omega}$, taking values $0$
and $1$. We require that the constraints be translation invariant, that
$c_x$ depend on a fixed finite neighborhood of $x$ and not on $\omega
_x$ [i.e., $c_x(\omega)=1$ if and only if $c_x(\omega^x)=1$]. We also
want the constraints to be monotone [if $\forall x\in\Z^d, \omega
_x\leq\omega'_x$, then $\forall x\in\Z^d, c_x(\omega)\geq
c_x(\omega
')$]. We will denote by $\L_{E}$ the generator of the environment
process: for $f$ a local function on $\lbrace0,1\rbrace^\Z$
%
%
\begin{equation}
\label{eq:generateurKCSM} \L_{E}f(\omega)=\sum_{y\in\Z}c_y(
\omega) \bigl((1-q) (1-\omega_y)+q\omega_y \bigr) \bigl[f
\bigl(\omega^y \bigr)-f(\omega) \bigr].
\end{equation}
In words, a zero (resp., each one) at site $x$ in configuration $\eta$
turns into a one (resp., a zero) at rate $(1-q)$ (resp., $q$), provided
the constraint is satisfied at $x$, that is, $c_x(\eta)=1$. This
process satisfies the detailed balance property w.r.t. $\mu$ the
product Bernoulli measure on $\lbrace0,1\rbrace^{\Z^d}$ of parameter
$1-q$, so it is reversible.

A transition $\omega\rightarrow\omega^x$ is \emph{legal} if
$c_x(\omega
)=1$. Note that $\omega\rightarrow\omega^x$ is legal if and only if
$\omega^x\rightarrow\omega$ is. A KCSM is \emph{noncooperative} if a
finite empty set is enough to empty the whole configuration through
legal transitions. More precisely, we have the following:

%
%
\begin{definition}\label{def:noncoop}
A KCSM is \emph{noncooperative} if the following holds:

There exists a finite set $A\subset\Z^d$ such that for every $\omega
\in
\Omega$, if $\omega_{|A}\equiv0$, for every $x\in\Z^d$ such that
$\omega_x=1$, there is a finite sequence $\omega^{(0)},\ldots,\omega
^{(n)}$ such that $\omega^{(0)}=\omega$, $ (\omega^{(n)}
)_x=0$, and for all $i=1,\ldots,n$, $\omega^{(i)}= (\omega
^{(i-1)} )^{x_i}$ where $x_i\in\Z^d$ such that $c_{x_i} (\omega
^{(i-1)} )=1$.
\end{definition}

The ergodic regime for KCSM was identified in \citet{KCSM}. In general,
there is a critical parameter $q_c\in[0,1]$ such that the process is
ergodic for $q
>q_c$ and nonergodic for $q<q_c$. $p_c=1-q_c$ is characterized as the
critical density of an appropriate bootstrap percolation model;
basically, it is the density above which blocked clusters (i.e.,
clusters of occupied sites that cannot be emptied through legal
transitions) appear with positive probability. A noncooperative model
is ergodic at every density $p=1-q\in(0,1)$ ($q_c=0$).

We now present the KCSM, which we will study in more detail.

\begin{figure}[b]

\includegraphics{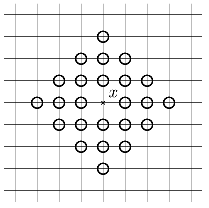}

\caption{$\mathcal{N}_3(x)$, the $3$-neighborhood of $x$ in $\Z^2$.}
\label{fig:3-neighborhood}
\end{figure}

We define a class of noncooperative KCSM, which we will call
``$k$-zeros'' for a positive integer $k$. Let $\|\cdot\|_1$ denote the
$1$-norm on $\Z^d$, that is, the norm induced by the graph distance. Let
%
%
\begin{equation}
\mathcal{N}_k(x)= \bigl\{y\in\Z^d | 0<\|y-x
\|_1\leq k \bigr\}
\end{equation}
be the $k$-neighborhood of $x$; see Figure~\ref{fig:3-neighborhood}.

%

The model ``$k$-zeros'' in $\Z^d$ is defined by the following
constraints [recall \eqref{eq:generateurKCSM}]:
%
%
\begin{equation}
\label{eq:contraintekzeros} c_x(\omega)= \cases{ %
1, &\quad $\mbox{if $ \displaystyle\sum_{y\in\mathcal{N}_k(x)}(1-
\omega_y)\geq k$},$
\vspace*{2pt}\cr
0, &\quad $\mbox{else},$}
\end{equation}
that is, the constraint is satisfied if there are at least $k$ zeros
within distance $k$. It is noncooperative since it is enough to empty
$0,e_1,2e_1,\ldots,(k-1)e_1$ to empty the whole lattice through legal
transitions. For $k=1$, the ``$1$-zero'' model is better known as the
one-flip Fredrickson--Andersen (or FA-1f) model.

The second model we want to study is the East model, a one-dimensional
KCSM for which the constraint is that the East neighbor of the
to-be-updated site be vacant. The corresponding generator is
%
%
\begin{equation}
\label{eq:eastgenerator} \mathcal{L}_Ef(\omega)=\sum
_{y\in\Z}(1-\omega_{y+1}) \bigl((1-q) (1-\omega
_y)+q\omega_y \bigr) \bigl[f \bigl(\omega^y
\bigr)-f(\omega) \bigr].
\end{equation}

In this study, we consider an environment given by a KCSM, and we
inject a tracer at its origin. The tracer jumps at rate one to each of
its nearest neighbors, provided that both the site where it sits and
the site where it wants to jump are empty (for the environment). More
formally, let $(\omega(t),X_t)$ be the joint evolution of the KCSM and
the tracer. It is a Markov process on $\lbrace0,1\rbrace^{\Z
^d}\times
\Z^d$ given by the generator
%
%
\begin{eqnarray}
\mathcal{L}_0f(\omega,x)&=&\sum_{y\in\Z^d}c_y(
\omega) \bigl((1-q) (1-\omega_y)+q\omega_y \bigr) \bigl[f
\bigl(\omega^y,x \bigr)-f(\omega,x) \bigr]
\nonumber
\\[-8pt]
\\[-8pt]
\nonumber
&&{} + \sum_{i=1}^d\sum
_{\alpha=\pm1}(1-\omega_x) (1-\omega_{x+\alpha e_i})
\bigl[f(\omega,x+\alpha e_i)-f(\omega,x) \bigr].
\end{eqnarray}

We consider the process $\eta(t)$ of the environment seen from the
tracer, whose generator is given by
%
%
\begin{eqnarray}\label{eq:processseenfromtracer}
\L f(\eta)&=&\sum_{y\in\Z^d}c_y(\eta)
\bigl((1-q) (1-\eta_y)+q\eta_y \bigr) \bigl[f \bigl(
\eta^y \bigr)-f(\eta) \bigr]
\nonumber
\\[-8pt]
\\[-8pt]
\nonumber
&&{} + \sum_{i=1}^d\sum
_{\alpha=\pm1}(1-\eta_0) (1-\eta_{\alpha
e_i})
\bigl[f(\eta_{\alpha e_i+\cdot})-f(\eta) \bigr],
\end{eqnarray}
where $\eta_{y+\cdot}$ denotes the configuration such that $ (\eta
_{y+\cdot} )_x=\eta_{y+x}$.
This is again a reversible process w.r.t. $\mu$ the product Bernoulli
measure on $\lbrace0,1\rbrace ^{\Z^d}$ of parameter $1-q$ (it
satisfies detailed balance).

A central tool in our study will be the spectral gap. Recall its definition.

%
\begin{definition}
The spectral gap of the generator $\mathcal{L}_E$ is given by the
variational principle
%
%
\begin{equation}
\label{eq:defgap} \gap(\mathcal{L}_E)=\inf\frac{-\mu(f\mathcal
{L}_Ef)}{\operatorname{Var}_\mu(f)},
\end{equation}
where the infimum is taken over all functions in $L^2(\mu)$ with
$\operatorname{Var}_\mu(f)\neq0$. A similar definition holds for $\gap(\mathcal{L})$
the spectral gap of the environment seen from the tracer.
\end{definition}

Recall also from \citet{aldousdiaconis,KCSM} that for the ``$k$-zeros''
model and the East model, the spectral gap is positive at any density.

\section{Main results}\label{sec:results}

We collect here the main results of this paper. The first one
establishes that after diffusive scaling the trajectory of the tracer
converges to a Brownian motion and introduces the diffusion coefficient
(or diffusion matrix) of the tracer.

%
\begin{prop}\label{prop:diffusion}
If the environment process is ergodic ($q>q_c$), we have
%
%
\begin{equation}
\lim_{\varepsilon\rightarrow0}\varepsilon X_{\varepsilon
^{-2}t}=\sqrt{2D}B_t,
\end{equation}
where $B_t$ is the standard Brownian motion, the convergence holds in
the sense of weak convergence of path measures on $D([0,\infty),\R^d)$
and the diffusion matrix $D$ is given by
%
%
\begin{equation}
\label{eq:formuleD} u.Du=q^2\|u\|_2^2-\int
_0^\infty\mu \bigl(j_{u}e^{\mathcal{L}t}j_{u}
\bigr)\,dt,
\end{equation}
where for any $u=(u_1,\ldots,u_d)\in\Z^d$ $j_u$ is given by the
action of
the generator $\mathcal{L}_0$ on the function $(\omega,x)\mapsto u.x$,
that is,
%
%
\begin{equation}
j_u(\eta)=(1-\eta_0)\sum_{i=1}^d
\sum_{\alpha=\pm1}(1-\eta_{\alpha
e_i})\alpha
u_i.
\end{equation}
\end{prop}

For the previous result to be meaningful, we need to prove $D>0$. In
the next proposition, we provide easy bounds on $D$ which show in
particular that this is true as soon as the KCSM has a positive
spectral gap. In \citet{KCSM}, it is proved for a large class of KCSM
that the spectral gap is positive in the whole ergodic regime, so this
requirement is not a big restriction. In particular, the spectral gap
is positive at every density $p=1-q\in(0,1)$ for the East model and
noncooperative models.

%
\begin{prop}\label{prop:easybounds}
%
%
\begin{equation}
q^2\|u\|_2^2\geq u.Du\geq
\frac{\gap(\mathcal{L}_E)}{4d+\gap(\mathcal
{L}_E)}q^2\|u\|_2^2.
\end{equation}
\end{prop}

The core of this paper is the study of $D$ when $q$ goes to zero, both
in noncooperative models and in the East model. In both cases the easy
bounds above can be significantly improved. For the sake of simplicity,
we give the following result only in the specific case of the
``$k$-zeros'' model. However, we expect our method to work more
generally for noncooperative models, and give the correct power of $q$
at high density.

%
\begin{theorem}\label{th:kzeros}
For the tracer diffusion in the ``$k$-zeros'' model, there exist
constants $0<c\leq C<\infty$ depending only on $d$ such that for all
$u\in\Z^d$,
%
%
\begin{equation}
cq^{k+1}\|u\|_2^2\leq u.Du\leq
Cq^{k+1}\|u\|_2^2.
\end{equation}
\end{theorem}

In the East model, we bound the ratio $D/\gap(\mathcal{L}_E
)$ on both sides by a polynomial in $q$.

%
\begin{theorem}\label{th:east}
When the environment is given by the East model, there exist constants
$C,c>0$ and $\alpha$ such that
%
%
\begin{equation}
\label{eq:east} cq^2\gap(\mathcal{L}_E)\leq D\leq
Cq^{-\alpha}\gap(\mathcal{L}_E).
\end{equation}
\end{theorem}

%
\begin{rem}
In~\citet{aldousdiaconis} and \citet{KCSM}, it is established that
%
%
\begin{equation}
\label{eq:gapeast} \lim_{q\rightarrow0}\frac{\log(1/\gap)}{ (\log
(1/q) )^2}=(2
\log2)^{-1}.
\end{equation}
In particular, this means that the powers of $q$ appearing in \eqref
{eq:east} are merely corrections to the correct asymptotic for $D$,
which is governed by the spectral gap of the East model. Inequality \eqref
{eq:east} is therefore incompatible with the prediction in \citet
{junggarrahanchandler} that $D\approx\gap^\xi$ for some $\xi<1$.
\end{rem}

\section{Convergence to a nondegenerate Brownian motion}\label
{sec:diffusivity}

We follow the strategy of \citet{kipnisvaradhan}, \citet
{demasiferrarigoldsteinwick} and \citet{spohn} to establish
Proposition~\ref{prop:diffusion}.

\begin{pf*}{Proof of Proposition~\ref{prop:diffusion}}
Considering the martingale
%
%
\begin{equation}
M^u_t=u.X_t-\int_0^t
j_u \bigl(\eta(s) \bigr)\,ds
\end{equation}
and following the steps of \citet{demasiferrarigoldsteinwick} and \citet
{spohn}, using reversibility, we get
%
%
\begin{eqnarray}\label{eq:limdiffusion}
&&\lim_{t\rightarrow\infty}\frac{1}{t}\mathbb{E}_{}
\bigl[{(u.X_t)^2}\bigr]
\nonumber
\\[-8pt]
\\[-8pt]
\nonumber
&&\qquad= \sum
_{i=1}^d\sum_{\alpha=\pm1}u_i^2
\mu \bigl((1- \eta_0) (1-\eta_{\alpha e_i}) \bigr)- 2\int
_0^\infty \mu \bigl(j_ue^{t\L}j_u
\bigr)\,dt.
\end{eqnarray}
In particular, $\int_0^\infty\mu(j_ue^{t\L}j_u )\,dt<\infty$,
so that, since the process of generator $\mathcal{L}$ is ergodic,
Theorem~1.8 of \citet{kipnisvaradhan} applies to $\int_0^t j_u(\eta
_s)\,ds$, yielding
%
%
\begin{equation}
\varepsilon u.X_{\varepsilon^{-2}t}=\varepsilon \bigl( M^u_{\varepsilon
^{-2}t}+N_{\varepsilon^{-2}t}
\bigr)+Q^\varepsilon(t),
\end{equation}
where $M_t+N_t$ is a martingale in $L^2(\mathbb{P})$ with stationary
increments, and $Q^\varepsilon(t)$ is an error term that vanishes when
$\varepsilon$ goes to $0$.
This implies the convergence of $\varepsilon X_{\varepsilon^{-2}t}$ to
$\sqrt{2D}B_t$ with $D$ given by \eqref{eq:formuleD}.
\end{pf*}

A first step in the direction of proving $D>0$ is to give a
variational formula for~$D$, which is the adaptation to our context of
Proposition~2 in \citet{spohn}.

%
\begin{lemme}\label{lem:varformula}
%
%
\begin{eqnarray}\quad
\label{eq:varformula} u.Du&=&\frac{1}{2}\inf_f \Biggl\{\sum
_{y\in\Z^d}\mu \bigl(c_y(\eta) \bigl((1-q) (1-
\eta_y)+q\eta_y \bigr) \bigl[f \bigl(\eta^y
\bigr)-f(\eta) \bigr]^2 \bigr)
\nonumber
\\[-8pt]
\\[-8pt]
\nonumber
&&\hspace*{28pt}{} + \sum_{i=1}^d\sum
_{\alpha
=\pm1} \mu \bigl((1-\eta_0) (1-
\eta_{\alpha e_i}) \bigl[\alpha u_i+f(\eta_{\alpha e_i+\cdot})-f(
\eta) \bigr]^2 \bigr) \Biggr\},
\end{eqnarray}
where the infimum is taken over local functions $f$ on $\Omega$.
\end{lemme}

\begin{pf}
We notice, as in \citet{spohn}, that
%
%
\begin{equation}
\label{eq:varformula1} \int_0^\infty\mu
\bigl(j_ue^{t\L}j_u \bigr)\,dt=-\inf \bigl\{-2
\mu(j_u f)-\mu(f\mathcal{L}f) \bigr\},
\end{equation}
where the infimum is taken over local functions on $\Omega$. Then,
using detailed balance, notice that we can write
%
%
\begin{equation}
\label{eq:var1} -4\mu(j_uf)= 2\sum_{i=1}^d
\sum_{\alpha=\pm1}\alpha u_i\mu \bigl((1-\eta
_0) (1-\eta_{\alpha e_i}) \bigl[f(\eta_{\alpha e_i+\cdot})-f(\eta)
\bigr] \bigr).
\end{equation}
Moreover,
%
%
\begin{eqnarray}\label{eq:var2}
-2\mu(f\mathcal{L}f)&=&\sum_{y\in\Z^d}\mu
\bigl(c_y(\eta) \bigl(p(1-\eta_y)+(1-p)\eta_y
\bigr) \bigl[f \bigl(\eta^y \bigr)-f(\eta) \bigr]^2 \bigr)
\nonumber
\\[-8pt]
\\[-8pt]
\nonumber
&&{} + \sum_{i=1}^d\sum
_{\alpha=\pm1} \mu \bigl((1-\eta_0) (1-\eta
_{\alpha e_i}) \bigl[f(\eta_{\alpha e_i+\cdot})-f(\eta) \bigr]^2
\bigr).
\end{eqnarray}
Inserting \eqref{eq:var1} and \eqref{eq:var2} into \eqref
{eq:limdiffusion} and rearranging the terms, we get \eqref{eq:varformula}.
\end{pf}

Now we can prove $D>0$ when the spectral gap of the environment is positive.

\begin{pf*}{Proof of Proposition~\ref{prop:easybounds}}
The upper bound follows directly from \eqref{eq:formuleD}, since the
second term is nonnegative.

For the lower bound, consider the expression of $D$ given in \eqref
{eq:varformula}. The first sum in the infimum is $-2\mu(f\mathcal
{L}_Ef)$, so that by definition of the spectral gap [recall~\eqref{eq:defgap}]
%
%
\begin{eqnarray}
\label{eq:lowerbound1}&& u.2Du\geq\inf \Biggl\{2\gap(\mathcal {L}_E)
\operatorname{Var}_\mu(f)
\nonumber
\\[-8pt]
\\[-8pt]
\nonumber
&&\hspace*{39pt}\qquad{}+
\sum_{i=1}^d\sum
_{\alpha=\pm1} \mu \bigl(\bar{\eta}_0\bar{
\eta}_{\alpha e_i} \bigl[\alpha u_i+f(\eta_{\alpha e_i+\cdot})-f(
\eta) \bigr]^2 \bigr) \Biggr\},
\end{eqnarray}
where we write $\bar{\eta}_x=1-\eta_x$.

To bound the double sum, we use the inequality $(a+b)^2\geq\gamma
a^2-\frac{\gamma}{1-\gamma}b^2$ for $\gamma<1$. This yields
\begin{eqnarray*}
&&\mu \bigl(\bar{\eta}_0\bar{\eta}_{\alpha e_i} \bigl[\alpha
u_i+f(\eta_{\alpha e_i})-f(\eta) \bigr]^2 \bigr)\\
&&\qquad
\geq \gamma q^2u_i^2-\frac
{\gamma}{1-\gamma}\mu
\bigl(\bar{\eta}_0\bar{\eta}_{\alpha e_i} \bigl[f(
\eta_{\alpha e_i+\cdot})-f(\eta) \bigr]^2 \bigr)
\\
&&\qquad\geq\gamma q^2 u_i^2-4
\frac{\gamma}{1-\gamma}\operatorname{Var}_\mu(f).
\end{eqnarray*}
So that, injecting this in \eqref{eq:lowerbound1}, we get
%
%
\begin{equation}
u.Du\geq\inf \biggl\{ \biggl(\gap(\mathcal{L}_E)-4d
\frac{\gamma}{1-\gamma
} \biggr)\operatorname{Var}_\mu(f)+\gamma q^2\|u
\|_2^2 \biggr\}.
\end{equation}
Choosing $\gamma=\frac{\gap(\mathcal{L}_E)}{4d+\gap(\mathcal{L}_E)}<1$,
we get the desired lower bound.
\end{pf*}

Note that at high density ($q\rightarrow0$), the spectral gap of the
East model is of order higher than any polynomial in $q$, so that the
term $q^2$ is negligible. In fact, for the East model, the lower bound
here is quite accurate (Theorem~\ref{th:east}). For noncooperative
models, however, we are able to do much better. In particular, for
FA-1f in one dimension, this gives $D\geq Cq^5$, which is pretty poor,
given that $D$ is in fact of order $q^2$, as predicted in \citet
{junggarrahanchandler}. Except in the FA-1f model, the upper bound also
needs refinement. Designing more precise bounds on $D$ when
$q\rightarrow0$ is the object of the next sections.

\section{Correct order of $D$ for small $q$ in noncooperative
models}\label{sec:kzeros}

%
\begin{rem}\label{rem:heuristic}
We believe that the techniques developed below can be adapted to show
the equivalent of Theorem~\ref{th:kzeros} for any noncooperative
model, $k$ being the minimal number of zeros needed to empty the whole
lattice (see Definition~\ref{def:noncoop}), and $1$ being replaced by
$m$ the minimal number of extra zeros needed to move a minimal cluster
around. We propose a heuristic for the order $q^{k+m}$, which we state
in dimension $1$ for simplicity. Consider for a moment a simple
symmetric random walk on the interval $\lbrace-1/(2q),\ldots,1/(2q)\rbrace
$ of length $1/q$. For large times $T$, the time spent in $0$ by the
random walk is approximately $Tq$. Since $1/q$ is the typical distance
between two zeros under the product Bernoulli measure $\mu$ on
$\lbrace
0,1\rbrace^\Z$, the fraction of time during which there is a zero at
$0$ before time $T$ is approximately~$Tq$. When that happens, a tracer
sitting in $0$ has a probability of order $q$ to jump, which gives a
diffusion coefficient for the tracer in the FA-1f model of order
$Tq\times q/T=q^2$. How does this adapt to another noncooperative
environment, where $k\geq1$, $m\geq1$ (e.g., the ``$k$-zeros''
model, $k> 1$, in which case $ m=1$)? A~single zero cannot move on its
own in such a model, but a group of $k$ zeros can, and since the number
of extra zeros it needs to move is $m$, the diffusion coefficient of
such a group is of order $q^m$. So we have to consider the fraction of
time spent in $0$ by a group of $k$ zeros performing a random walk on
$\lbrace-1/(2q^k),\ldots,1/(2q^k)\rbrace$ before time $T$ ($1/q^k$ being
the typical distance between two such groups under $\mu$), that is,
$Tq^k$. During the time the group of $k$ zeros is in contact with the
tracer (i.e., at site $0$), the tracer diffuses with it, which means
with rate $q^m$. In the end, the diffusion coefficient of the tracer
should therefore be of order $Tq^k\times q^m/T$.
\end{rem}

\subsection{Lower bound in Theorem~\texorpdfstring{\protect\ref{th:kzeros}}{3.3}}\label{sec:lowerbound}

The key to the proof of the lower bound we give below is that we are
able to come down to studying a \emph{local} dynamics; see Lemma~\ref
{lem:Dbarre} and the description of the dynamics in the proof of
Lemma~\ref{lem:lowerboundDbarre}. The possibility of doing this
simplification is strongly related to the fact that we are working with
noncooperative models.

For the sake of simplicity, this proof is written for $k=3$, but it
generalizes without difficulty to any $k\geq1$. It is widely inspired
by the fourth section in \citet{spohn}.

\begin{figure}

\includegraphics{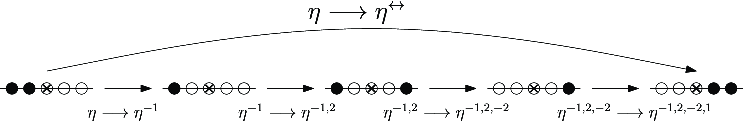}

\caption{The four legal flips used to reconstruct the swap $\eta
\longrightarrow\eta^\leftrightarrow$ when
$\eta_{-1}=\eta_{-2}=\bar{\eta}_1= (1-\eta_2)=1$. The cross recalls that
the tracer is sitting at the origin.}
\label{fig:swapinflips}
\end{figure}

\begin{figure}[b]

\includegraphics{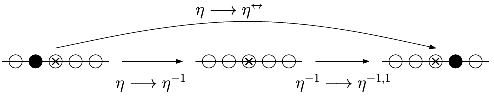}

\caption{The two legal flips used to reconstruct the swap $\eta
\longrightarrow\eta^\leftrightarrow$ when
$\eta_{-1}=\bar{\eta}_{-2}=\bar{\eta}_1= \bar{\eta}_2=1$.}
\label{fig:swapinflips2}
\end{figure}

The first step is to give a lower bound on $D$ in terms of the
diffusion coefficient $\overline{D}$ of another dynamics (Lemma~\ref
{lem:Dbarre}), for which we can prove positivity (Lemma~\ref
{lem:lowerboundDbarre}). In the auxiliary dynamics, the only allowed
transitions are jumps of the tracer between empty sites and swaps of
its left and right neighborhood, which can be reconstructed using only
transitions that are allowed in the initial dynamics; see Figures~\ref
{fig:swapinflips} and \ref{fig:swapinflips2}. We need some notation to
be more specific.

%

%

Let $\mu^{(3)}$ be the product Bernoulli measure \emph{on $\Z$}
conditioned to having at least three consecutive zeros, one of which at
the origin; that is, let $A\subset\Omega$ be defined as
%
%
\begin{eqnarray}\quad
\label{eq:defA} &&A= \bigl\{\eta\in\Omega| \eta_0=0
\nonumber
\\[-8pt]
\\[-8pt]
\nonumber
&&\hspace*{5pt}\qquad\mbox{and }(1-
\eta_1) (1-\eta_2)+(1-\eta_{-1}) (1-
\eta_1)+(1-\eta_{-2}) (1-\eta_{-1})\geq1 \bigr\}
\end{eqnarray}
and
%
%
\begin{equation}
\label{eq:defmu3} \mu^{(3)}=\mu (\cdot| A ).
\end{equation}
Also, if $\eta\in\Omega$, denote by $\eta^{\leftrightarrow}$ the
configuration obtained by exchanging the occupation numbers in sites
$-1$ and $+1$, and $-2$ and $+2$
%
%
\begin{equation}
\label{eq:defexchange} \eta^{\leftrightarrow}_y=
\cases{\eta_1, &\quad $\mbox{if $y=-1$},$
\vspace*{2pt}\cr
\eta_{-1}, &\quad $\mbox{if $y=1$},$
\vspace*{2pt}\cr
\eta_2, &\quad $\mbox{if $y=-2$},$
\vspace*{2pt}\cr
\eta_{-2}, &\quad $\mbox{if $y=2$},$
\vspace*{2pt}\cr
\eta_y, &\quad  $\mbox{else.}$}
\end{equation}
We also generalize the notation $\eta^x$ by defining $\eta
^{x_1,\ldots,x_n}$ as the configuration $\eta$ flipped at sites
$x_1,\ldots,x_n$ (the $x_i$ being distinct).

We can now state the following:

%
%
\begin{lemme}\label{lem:Dbarre}
If $\overline{D}$ is defined by
%
%
\begin{eqnarray}
\label{eq:defDbarre} \overline{D}&=&\frac{1}{2}\inf_f
\bigl\{\mu^{(3)} \bigl( \bigl(1-(1-\eta_1) (1-
\eta_{-1}) \bigr) \bigl[f \bigl(\eta^\leftrightarrow \bigr)-f(\eta)
\bigr]^2 \bigr)
\nonumber
\\
&&\hspace*{24pt}{}+ \mu^{(3)} \bigl((1-\eta_1) \bigl[1+f(
\eta_{1+\cdot})-f(\eta) \bigr]^2 \bigr)\\
&&\hspace*{38pt}{}+\mu^{(3)}
\bigl((1-\eta_{-1}) \bigl[-1+f(\eta_{-1+\cdot
})-f(\eta)
\bigr]^2 \bigr) \bigr\},\nonumber
\end{eqnarray}
where the infimum is taken over local functions on $\Omega$, then we have
%
%
\begin{equation}
e_1.De_1\geq\frac{1+2p}{4}q^4
\overline{D}.
\end{equation}
\end{lemme}

\begin{pf}
For briefness, we define
%
%
\begin{equation}
\label{eq:defrate} \bar{\eta}_x= 1-\eta_x \quad\mbox{and}\quad
r_x(\eta)=(1-q)\bar{\eta}_x+q\eta_x.
\end{equation}
Then we have, given the definition of $\mu^{(3)}$ \eqref{eq:defmu3},
for every local function $f$,
%
%
\begin{eqnarray}\label{eq:lemDbarre1}
&&\mu^{(3)} \bigl( (1-\bar{\eta}_1\bar{\eta}_{-1}
) \bigl[f \bigl(\eta^\leftrightarrow \bigr)-f(\eta) \bigr]^2 \bigr)
\nonumber\\
&&\qquad= \mu^{(3)} \bigl(\bar{\eta}_1\bar{\eta}_2
\eta_{-1} \bigl[f \bigl(\eta^\leftrightarrow \bigr)-f(\eta)
\bigr]^2 \bigr)
\\
&&\qquad\quad {}+ \mu^{(3)} \bigl(\bar{\eta}_{-1}\bar{\eta}_{-2}
\eta_1 \bigl[f \bigl(\eta^\leftrightarrow \bigr)-f(\eta)
\bigr]^2 \bigr).\nonumber
\end{eqnarray}
Our aim is to reconstruct the swap changing $\eta$ into $\eta
^\leftrightarrow$, using only legal (for the ``$3$-zeros'' model
dynamics) flips. The first term of the RHS in \eqref{eq:lemDbarre1} can
be rewritten as
%
%
\begin{eqnarray}
\label{eq:lemDbarre3} &&\mu^{(3)} \bigl(\bar{\eta}_1\bar{
\eta}_2\eta_{-1}\eta_{-2} \bigl[f \bigl(
\eta^\leftrightarrow \bigr)-f(\eta) \bigr]^2 \bigr)
\nonumber
\\[-8pt]
\\[-8pt]
\nonumber
&&\qquad{}+\mu^{(3)}
\bigl(\bar{\eta}_1\bar{\eta}_2\eta_{-1}\bar{
\eta}_{-2} \bigl[f \bigl(\eta^\leftrightarrow \bigr)-f(\eta)
\bigr]^2 \bigr).
\end{eqnarray}
Let us focus on the first term. See in Figure~\ref{fig:swapinflips} a
representation of the successive flips used to reconstruct the swap.
Writing that, when $\eta_{-1}=\eta_{-2}=\bar{\eta}_1=\bar{\eta}_2=1$,
%
%
\begin{eqnarray}\qquad
f \bigl(\eta^\leftrightarrow \bigr)-f(\eta)&=&f \bigl(\eta^{-1,2,-2,1}
\bigr)-f \bigl(\eta^{-1,2,-2} \bigr)+f \bigl(\eta^{-1,2,-2} \bigr)-f
\bigl(\eta^{-1,2} \bigr)
\nonumber
\\[-8pt]
\\[-8pt]
\nonumber
&&{}+ f \bigl(\eta^{-1,2} \bigr)-f \bigl(\eta^{-1} \bigr)+f
\bigl(\eta^{-1} \bigr)-f (\eta),
\end{eqnarray}
and using the Cauchy--Schwarz inequality, we have
%
%
\begin{eqnarray}
&&\mu^{(3)} \bigl(\bar{\eta}_1\bar{\eta}_2
\eta_{-1}\eta_{-2} \bigl[f \bigl(\eta^\leftrightarrow \bigr)-f(
\eta) \bigr]^2 \bigr)\nonumber\\
&&\qquad\leq 4\mu^{(3)} \bigl(\bar{
\eta}_1\bar{\eta}_2\eta_{-1}\eta_{-2}
\bigl[f \bigl(\eta^{-1,2,-2,1} \bigr)-f \bigl(\eta^{-1,2,-2} \bigr)
\bigr]^2 \bigr)
\nonumber
\\
&&\qquad\quad{}+ 4\mu^{(3)} \bigl(\bar{\eta}_1\bar{\eta}_2
\eta_{-1}\eta_{-2} \bigl[f \bigl(\eta^{-1,2,-2} \bigr)-f
\bigl(\eta^{-1,2} \bigr) \bigr]^2 \bigr)
\\
&&\qquad\quad{}+ 4\mu^{(3)} \bigl(\bar{\eta}_1\bar{\eta}_2
\eta_{-1}\eta_{-2} \bigl[f \bigl(\eta^{-1,2} \bigr)-f
\bigl(\eta^{-1} \bigr) \bigr]^2 \bigr)\nonumber
\\
&&\qquad\quad{}+ 4\mu^{(3)} \bigl(\bar{\eta}_1\bar{\eta}_2
\eta_{-1}\eta_{-2} \bigl[f \bigl(\eta^{-1} \bigr)-f
(\eta) \bigr]^2 \bigr).\nonumber
\end{eqnarray}
Note that all the flips involved are legal for the dynamics
``$3$-zeros'': there are always at least three zeros in the
$3$-neighborhood of the site that is flipped. Then we make a change of
variables in the first three terms above to get
%
%
\begin{eqnarray}\label{eq:lemDbarre2}
&&\mu^{(3)} \bigl(\bar{\eta}_1\bar{\eta}_2
\eta_{-1}\eta_{-2} \bigl[f \bigl(\eta^\leftrightarrow \bigr)-f(
\eta) \bigr]^2 \bigr)\nonumber\\
&&\qquad\leq 4\frac
{1-q}{q}\mu^{(3)}
\bigl(\bar{\eta}_1\eta_2\bar{\eta}_{-1}\bar{
\eta}_{-2} \bigl[f \bigl(\eta^{1} \bigr)-f (\eta)
\bigr]^2 \bigr)
\nonumber
\\
&&\qquad\quad{}+ 4\mu^{(3)} \bigl(\bar{\eta}_1\eta_2\bar{
\eta}_{-1}\eta_{-2} \bigl[f \bigl(\eta^{-2} \bigr)-f
(\eta) \bigr]^2 \bigr)
\\
&&\qquad\quad{}+ 4\frac{1-q}{q}\mu^{(3)} \bigl(\bar{\eta}_1\bar{
\eta}_2\bar{\eta}_{-1}\eta_{-2} \bigl[f \bigl(
\eta^{2} \bigr)-f (\eta) \bigr]^2 \bigr)\nonumber
\\
&&\qquad\quad{}+ 4\mu^{(3)} \bigl(\bar{\eta}_1\bar{\eta}_2
\eta_{-1}\eta_{-2} \bigl[f \bigl(\eta^{-1} \bigr)-f
(\eta) \bigr]^2 \bigr).\nonumber
\end{eqnarray}

In the same way (following the strategy represented in Figure~\ref
{fig:swapinflips2}), we get
%
%
\begin{eqnarray}
\label{eq:lemDbarre5}&& \mu^{(3)} \bigl(\bar{\eta}_1\bar{
\eta}_2\eta_{-1}\bar{\eta}_{-2} \bigl[f \bigl(
\eta^\leftrightarrow \bigr)-f(\eta) \bigr]^2 \bigr)\nonumber\\
&&\qquad\leq
\frac
{2(1-q)}{q}\mu^{(3)} \bigl(\bar{\eta}_{-1}\bar{
\eta}_2\bar{\eta}_{-2}\bar{\eta}_1 \bigl[f
\bigl(\eta^1 \bigr)-f(\eta) \bigr]^2 \bigr)
\\
&& \qquad\quad{}+ {2}\mu^{(3)} \bigl(\bar{\eta}_1\bar{
\eta}_2\bar{\eta}_{-2}\eta_{-1} \bigl[f \bigl(
\eta^{-1} \bigr)-f(\eta) \bigr]^2 \bigr).\nonumber
\end{eqnarray}

Combining \eqref{eq:lemDbarre3}, \eqref{eq:lemDbarre2} and \eqref
{eq:lemDbarre5} and doing the same for the second term in \eqref
{eq:lemDbarre1}, we get [recall \eqref{eq:defrate} for the definition
of $r_x$]
%
%
\begin{eqnarray}\label{eq:lemDbarre4}
&&\mu^{(3)} \bigl( (1-\bar{\eta}_1\bar{\eta}_{-1}
) \bigl[f \bigl(\eta^\leftrightarrow \bigr)-f(\eta) \bigr]^2 \bigr)\nonumber\\
&&\qquad
\leq \frac
{4}{q}\mu^{(3)} \bigl(\bar{\eta}_{-1}\bar{
\eta}_{-2}r_1(\eta) \bigl[f \bigl(\eta^{1}
\bigr)-f (\eta) \bigr]^2 \bigr)
\\
&&\qquad\quad{}+ \frac{4}{q}\mu^{(3)} \bigl(\bar{\eta}_1\bar{
\eta}_{-1}r_{-2}(\eta) \bigl[f \bigl(\eta^{-2}
\bigr)-f (\eta) \bigr]^2 \bigr)\nonumber
\\
&&\qquad\quad{}+ \frac{4}{q}\mu^{(3)} \bigl(\bar{\eta}_1\bar{
\eta}_{-1}r_2(\eta) \bigl[f \bigl(\eta^{2}
\bigr)-f (\eta) \bigr]^2 \bigr)
\nonumber
\\
&&\qquad\quad{}+ \frac{4}{q}\mu^{(3)} \bigl(\bar{\eta}_1\bar{
\eta}_2r_{-1}(\eta) \bigl[f \bigl(\eta^{-1}
\bigr)-f (\eta) \bigr]^2 \bigr).\nonumber
\end{eqnarray}

Now notice that we have
\[
\mu^{(3)} \bigl(\bar{\eta}_{-1}\bar{\eta}_{-2}r_1(
\eta) \bigl[f \bigl(\eta^{1} \bigr)-f (\eta) \bigr]^2
\bigr) = \frac{1}{\mu(A)}\mu \bigl(\bar{\eta}_0\bar{
\eta}_{-1}\bar{\eta}_{-2}r_1(\eta) \bigl[f \bigl(
\eta^{1} \bigr)-f (\eta) \bigr]^2 \bigr)
\]
and similarly for the other terms in \eqref{eq:lemDbarre4}, so that we
have proved the following inequality, recalling that $\mu(A)=q^3(1+2p)$:
%
%
\begin{eqnarray}
\label{eq:finlemDbarre1}&& \sum_{y\in\Z^d}\mu \bigl(c_y(
\eta)r_y(\eta) \bigl[f \bigl(\eta^y \bigr)-f(\eta)
\bigr]^2 \bigr)
\nonumber
\\[-8pt]
\\[-8pt]
\nonumber
&&\qquad\geq q^4\frac{(1+2(1-q))}{4}
\mu^{(3)} \bigl( (1-\bar{\eta}_1\bar{\eta}_{-1} )
\bigl[f \bigl(\eta^\leftrightarrow \bigr)-f(\eta) \bigr]^2 \bigr).
\end{eqnarray}

We are almost done. It remains to notice that
\begin{eqnarray}
\mu^{(3)} \bigl(\bar{\eta}_1 \bigl[1+f(
\eta_{1+\cdot})-f(\eta) \bigr]^2 \bigr)&\leq&\frac{1}{\mu(A)}
\mu \bigl(\bar{\eta}_0\bar{\eta}_1 \bigl[1+f(
\eta_{1+\cdot})-f(\eta) \bigr]^2 \bigr)
\nonumber
\end{eqnarray}
and similarly with $1$ replaced by $-1$, so that a fortiori
\begin{eqnarray}\label{eq:finlemDbarre2}
&&\mu^{(3)} \bigl(\bar{\eta}_1 \bigl[1+f(\eta_{1+\cdot})-f(\eta)
\bigr]^2 \bigr)+\mu^{(3)} \bigl(\bar{\eta}_{-1} \bigl[-1+f(\eta_{-1+\cdot
})-f(\eta) \bigr]^2 \bigr)
\nonumber
\\[-8pt]
\\[-8pt]
\nonumber
 &&\qquad\leq\frac{4}{q^4(1+2(1-q))}\sum_{i=1}^d
\sum_{\alpha=\pm1} \mu \bigl(\bar{\eta}_0\bar{
\eta}_{\alpha e_i} \bigl[\alpha\delta_{1i}+f(\eta
_{\alpha e_i+\cdot})-f(\eta) \bigr]^2 \bigr).
\end{eqnarray}

Combining \eqref{eq:finlemDbarre1} and \eqref{eq:finlemDbarre2}, and
recalling \eqref{eq:varformula}, we get the lemma.
\end{pf}

Of course there is nothing special about the direction $e_1$, and the
lemma is valid in all directions. Notice that it does not depend on the
dimension. We now complete the proof of the lower bound in Theorem~\ref
{th:kzeros} by providing a universal lower bound on $\overline{D}$.

%
\begin{lemme}\label{lem:lowerboundDbarre}
$\overline{D}$ defined in \eqref{eq:defDbarre} is the diffusion
coefficient of a universal auxiliary dynamics and is bounded below as
%
%
\begin{equation}
\overline{D}\geq4/9.
\end{equation}
\end{lemme}

\begin{pf}
Following the same lines as in the proof of Proposition~\ref
{prop:diffusion} and Lemma~\ref{lem:varformula}, we see that
$\overline
{D}$ is the diffusion coefficient of the dynamics reversible w.r.t.
$\mu
^{(3)}$ described below:
\begin{itemize}
\item with rate $1$, if $\eta_1=0$, the tracer jumps to the right, that
is, we go from $\eta$ to $\eta_{1+\cdot}$,
\item with rate $1$, if $\eta_{-1}=0$, the tracer jumps to the left,
that is, we go from $\eta$ to $\eta_{-1+\cdot}$,
\item with rate $1$, if either $\eta_1=1$ or $\eta_{-1}=1$, $\lbrace
-2,-1\rbrace$ and $\lbrace2,1\rbrace$ are swapped, that is, we go from
$\eta$ to $\eta^{\leftrightarrow}$.
\end{itemize}

As in \citet{spohn}, starting from a configuration $\eta$ chosen in $A$
[recall \eqref{eq:defA}], we can index by $\Z$ all the configurations
that can be reached by this dynamics in the following way. $\eta
^{(0)}=\eta$ is the initial configuration, that is almost surely in
$A$. Then we define inductively $\eta^{(n)}$, $n\in\Z$. If $\eta
^{(n)}_1=0$, $\eta^{(n+1)}=\eta^{(n)}_{1+\cdot}$. If $\eta^{(n)}_1=1$,
$\eta^{(n+1)}= (\eta^{(n)} )^\leftrightarrow$. Similarly, if
$\eta^{(n)}_{-1}=0$, $\eta^{(n-1)}=\eta^{(n)}_{-1+\cdot}$. If $\eta
^{(n)}_{-1}=1$, $\eta^{(n-1)}= (\eta^{(n)} )^\leftrightarrow$.
Note that this definition is consistent ($\eta^{(n+1-1)}=\eta^{(n)}$).

Using this labeling with integers of all attainable configurations, the
dynamics described above can be equivalently defined in the following
way: if the system is in the configuration $\eta^{(n)}$, it goes to
$\eta^{(n+1)}$ with rate one, and to $\eta^{(n-1)}$ also with rate one.
So we can rewrite the process starting from $\eta$ as $\eta(t)=\eta
^{(N_t)}$ where $(N_t)_{t\geq0}$ is a simple random walk on $\Z$.

Now to conclude, we just need to notice that if $X_t$ is the position
of the tracer at time $t$ in this dynamics, we have
\[
|X_t|\geq \bigl\lfloor\tfrac{2}{3}|N_t| \bigr
\rfloor,
\]
since two out of three times $N$ moves to the right, $X$ also jumps by
one (and similarly to the left).
\[
2\overline{D}=\lim_{t\rightarrow+\infty}\frac{1}{t}\mathbb
{E}_{}\bigl[{X_t^2}\bigr]\geq
\frac
{4}{9}\lim_{t\rightarrow+\infty} \frac{1}{t}
\mathbb{E}_{}\bigl[{N_t^2}\bigr]=8/9.
\]
\upqed\end{pf}

To deduce Theorem~\ref{th:kzeros} lower bound from Lemma~\ref
{lem:Dbarre} and Lemma~\ref{lem:lowerboundDbarre}, let $u\in\R^d$ be
such that $\|u\|_2=1$ and notice that we can use comparisons with the
auxiliary dynamics above in all directions to get
%
%
\begin{eqnarray}
2u.Du & \geq& \sum_{i=1}^d\inf
_{f_i} \biggl\{\frac{1}{d}\sum
_{x\in\Z
^d}\mu \bigl(c_x(\eta)r_x(\eta)
\bigl[f_i \bigl(\eta^x \bigr)-f_i(\eta)
\bigr]^2 \bigr)
\nonumber
\\
&&\hspace*{33pt}{} + \sum_{\alpha=\pm
1}\mu \bigl((1- \eta_0)
(1-\eta_{\alpha e_i}) \bigl[\alpha u_i+f_i(
\eta_{\alpha e_i+\cdot})-f_i(\eta) \bigr]^2 \bigr) \biggr\}
\nonumber
\\
&\geq& \sum_{i=1}^d u_i^2
\inf_{f_i} \biggl\{\frac{1}{d}\sum
_{x\in\Z\cdot
e_i}\mu \bigl(c^{i}_x(
\eta)r_x(\eta) \bigl[f_i \bigl(\eta^x
\bigr)-f_i(\eta) \bigr]^2 \bigr)
\\
&&\hspace*{44pt}{} + \frac{1}{d}\sum_{\alpha=\pm
1}\mu \bigl((1-
\eta_0) (1-\eta_{\alpha e_i}) \bigl[\alpha+f_i(
\eta_{\alpha
e_i+\cdot})-f_i(\eta) \bigr]^2 \bigr) \biggr\}
\nonumber
\\
&\geq& \frac{2}{d} D_1,\nonumber
\end{eqnarray}
where $c^i_x(\eta)$ is one if and only if the constraint is satisfied
using only zeros in the direction $i$, $D_1$ is the diffusion
coefficient in one dimension and we used $\sum_{i=1}^d u_i^2=1$.
Theorem~\ref{th:kzeros} follows from this inequality and the two
previous lemmas.

%
\begin{rem}
This strategy can be applied to other noncooperative models. However,
the auxiliary dynamics (the one involving swaps around the origin and
jumps of the tracer) will be model dependant and may not be strictly
one-dimensional. It may be encoded by a random walks on graphs slightly
more complex than $\Z$, but still with a uniformly positive diffusion
coefficient. We believe that this technique could allow us to retrieve
the correct exponent at low temperature for noncooperative models.
\end{rem}

\subsection{Upper bound in Theorem~\texorpdfstring{\protect\ref{th:kzeros}}{3.3}}

In view of~\eqref{eq:varformula}, to find an upper bound on~$D$, we
need to find an appropriate test function. As a warming, suppose that
$d=1$. Then, looking for a function that cancels the second line
in~\eqref{eq:varformula}, we find that a natural function to consider is
%
%
\begin{equation}
\label{eq:deffonctiontestdim1} f(\eta)=\min\{x\in\N| \eta_x=1 \}.
\end{equation}
Then it is not too difficult to check that if we plug this function
into the first line of~\eqref{eq:varformula}, we get an expression of
order $q^{k+1}$: the factor $q^k$ comes from the constraint, and the
extra $q$ comes from the extra empty site we need in order to evolve.

\begin{figure}[b]

\includegraphics{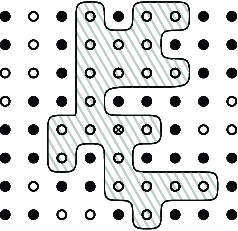}

\caption{An example of $C(\eta)$. Zeros are represented by empty
circles, ones by filled
disks and the origin is marked by a cross. The cluster of zeros
containing the origin is
circled by a line and tiled in gray. In this case, $f(\eta)=4$.}
\label{fig:fonctiontest}
\end{figure}

In higher dimension, we are going to find a good test function to
evaluate $e_1.De_1$. Define $C(\eta)$ the connected cluster of zeros
containing the origin in the configuration $\eta$ [$C(\eta
)=\varnothing
$ if $\eta_0=1$]. See Figure~\ref{fig:fonctiontest} for an example.

%

Now we can define our test function.
%
%
\begin{equation}
f(\eta)=\min \bigl\{x\in\N| C(\eta)\subset(-\infty,x-1 ]\times \Z^{d-1}
\bigr\}.
\end{equation}
For instance, if $\eta_0=1$, $f(\eta)=0$. In Figure~\ref
{fig:fonctiontest}, $f(\eta)=4$. Note that this function coincides with
that in \eqref{eq:deffonctiontestdim1} when $d=1$. This function
cancels the second line in~\eqref{eq:varformula} when $u=e_1$. Indeed,
when $(1-\eta_0)(1-\eta_{\alpha e_i})\neq0$, $0$ and $\alpha e_i$
belong to the same cluster of zeros. So what we need to do is show that
%
%
\begin{equation}
\label{eq:boundflipterms} \sum_{y\in\Z^d}\mu \bigl(c_y(
\eta) \bigl((1-q) (1-\eta_y)+q\eta_y \bigr) \bigl[f
\bigl( \eta^y \bigr)-f(\eta) \bigr]^2 \bigr)\leq
Cq^{k+1}
\end{equation}
for some finite $C$. Let us split the LHS into two terms and treat them
separately: we need to show that
%
%
\begin{eqnarray}
S_0& = & \sum_{y\in\Z^d}\mu
\bigl(c_y(\eta) (1-\eta_y) \bigl[f \bigl(\eta
^y \bigr)-f(\eta) \bigr]^2 \bigr)\leq Cq^{k+1},
\label{eq:boundS0}
\\
S_1 &= & \sum_{y\in\Z^d}\mu
\bigl(c_y(\eta)\eta_y \bigl[f \bigl(\eta^y
\bigr)-f(\eta) \bigr]^2 \bigr)\leq Cq^k\label{eq:boundS1}.
\end{eqnarray}
Thanks to detailed balance, $(1-q)S_0=qS_1$, so we only need to show
\eqref{eq:boundS0}.

Let us now study $S_0$. The mechanism involved here is the removal of
part of the cluster of zeros around the origin. In particular, when
$(1-\eta_y) [f(\eta^y)-f(\eta) ]^2\neq0$, we certainly have
\[
\bigl[f \bigl(\eta^y \bigr)-f(\eta) \bigr]^2\leq \bigl
\llvert C(\eta) \bigr\rrvert^2,
\]
where $\llvert C(\eta)\rrvert$ is the cardinal of $C(\eta)$. So that
%
%
\begin{eqnarray}
S_0&\leq& \mu \biggl( \bigl\llvert C(\eta) \bigr\rrvert
^2\sum_{y\in C(\eta)}c_y(\eta) (1-
\eta_y) \biggr)
\nonumber
\\[-8pt]
\\[-8pt]
\nonumber
&\leq& \sum_{n\geq0}\mu \biggl( \bigl\llvert C(\eta)
\bigr\rrvert^2\mathbf{1}_{0\leftrightarrow\partial B_n,
0\nleftrightarrow\partial
B_{n+1}}\sum
_{y\in C(\eta)}c_y(\eta) (1-\eta_y) \biggr),
\end{eqnarray}
where $\partial B_n$ denotes the set of points at distance $n$ from
$0$, and $ \lbrace0\leftrightarrow\partial B_n \rbrace$ is
the event that there is a site at distance $n$ from $0$ in $C(\eta)$.
Since on the event $\lbrace0\leftrightarrow\partial B_n,
0\nleftrightarrow\partial B_{n+1}\rbrace$, $C(\eta)\subset B_1(0,n)$,
we have
%
%
\begin{equation}
S_0  \leq \sum_{n\geq0}(2n+1)^{2d}
\sum_{y\in B_1(0,n)}\mu \bigl(c_y(\eta) (1-
\eta_y)\mathbf{1}_{0\leftrightarrow\partial B_n, 0\nleftrightarrow
\partial B_{n+1}} \bigr).
\end{equation}
On the one hand, for any $y$, we have for some constant $C$ depending
only on $d$,
%
%
\begin{equation}
\label{eq:upper1} \mu \bigl(c_y(\eta) (1-\eta_y) \bigr)\leq
Cq^{k+1},
\end{equation}
since the constraint requires at least $k$ zeros to be satisfied, and
$c_y$ is independent from $\eta_y$. On the other hand, if
$0\leftrightarrow\partial B_n$, there is a self-avoiding walk of
length $n$ starting at $0$ which is empty. So a rough bound on the
number of self-avoiding walks of length $n$ yields
%
%
\begin{equation}
\label{eq:upper2} \mu({0\leftrightarrow\partial B_n, 0\nleftrightarrow
\partial B_{n+1}} )\leq(2d)^n q^n.
\end{equation}
Putting together \eqref{eq:upper1} and \eqref{eq:upper2}, we get
%
%
\begin{equation}
S_0\leq \sum_{n\geq0}(2n+1)^{3d}
\bigl[Cq^{k+1}\wedge(2dq)^n \bigr]\leq C'q^{k+1}
\end{equation}
for $q$ small enough. So we have proved \eqref{eq:boundS0}.

A general argument allows to retrieve the upper bound in Theorem~\ref
{th:kzeros} for any $u\in\R^d$ from the result for $e_1,\ldots,e_d$. Write
$u=\sum_{i=1}^du_ie_i$, and compute
%
%
\begin{equation}
u.Du=\sum_{i=1}^d u_i^2
e_i.De_i+\sum_{i\neq j}u_iu_je_i.De_j.
\end{equation}
Notice that $D$ is symmetric and positive (by Proposition~\ref
{prop:easybounds}), so that the application $(u,v)\mapsto u.Dv$ is a
scalar product. We can therefore apply the Cauchy--Schwarz inequality
to the terms $e_i.De_j$ and get
%
%
\begin{equation}
u.Du\leq Cq^{k+1} \Biggl(\sum_{i=1}^d|u_i|
\Biggr)^2 \leq C' q^{k+1},
\end{equation}
where $C'$ depends only on $d$ by equivalence of the norms in finite dimension.

\section{In the East model, \texorpdfstring{$D\approx\gap$}{D approx gap}}\label{sec:east}

In this section, we prove Theorem~\ref{th:east}.

Before getting into the results concerning the tracer, let us recall
briefly the definition and basic property of the so-called
distinguished zero, a very useful tool for the study of the East model,
which was introduced in \citet{aldousdiaconis}.

%
\begin{definition}
Consider $\omega\in\Omega$ a configuration with $\omega_x=0$ for some
\mbox{$x\in\mathbb{Z}$}. Define $\xi(0)=x$. Call $T_1=\inf\lbrace t\geq0
| \mbox{\emph{the clock in }}x\mbox{\emph{ rings and }}\omega
_{x+1}(t)=0\rbrace$, the time of the first legal ring at $x$. Let $\xi
(s)=x$ for $s<T_1$, $\xi(T_1)=x+1$ and start again to define
recursively $ (\xi(s) )_{s\geq0}$.
\end{definition}

Notice that for any $s\geq0$, $\omega_{\xi(s)}(s)=0$, and that $\xi
\dvtx
\mathbb{R}^+\rightarrow\Z$ is almost surely c\`adl\`ag and increasing
by jumps of $1$.

This distinguished zero has an important property: as it moves forward,
it leaves equilibrium on its left; see Lemma~4 in \citet{aldousdiaconis}
or Lemma~3.5 of \citet{nonequilibrium}. In particular, if $\omega$ is
such that $\omega_x=0$ and $A$ an event depending only on the
configuration restricted to $[x_-,x_+]$, with $x_+<x$, letting
$V=\lbrace x_-,\ldots,x-1\rbrace$, then we have the following estimate:
%
%
\begin{equation}
\label{eq:distinguishedestimate} \mathbb{P}_{\omega}\bigl({\omega(t)\in A}\bigr)\leq
\mu_V( \omega_{|V})^{-1}\mathbb{P}_{\mu_V\cdot\omega}
\bigl({\omega(t)\in A}\bigr)=\mu(\omega_{|V})^{-1}\mu(A),
\end{equation}
where $\mu_{V}$ is the Bernoulli($1-q$) product measure on $\lbrace
0,1\rbrace^V$, $\mu_{V}\cdot\omega$ denotes the law of a random
configuration equal to $\omega$ on $\Z\setminus V$ and chosen with law
$\mu_{V}$ on $V$. In the above estimate, the factor $\mu_V(\omega
_{|V})$ comes from a change of measure to start from $\mu$ in $V$, and
the last equality comes from the property of the distinguished zero
mentioned above.

For briefness, in this section, we will denote the spectral gap of the
East process by $\gap$; see~(\ref{eq:defgap}).

\begin{pf*}{Proof of Theorem~\ref{th:east}}
The lower bound is already contained in Proposition~\ref{prop:easybounds}.

For the proof of the upper bound, fix $t>0$ and $\tau\ll t$ to be
chosen later, such that $t/\tau$ is an integer and $\tau\lesssim\gap
^{-1}$ (more precisely, $\tau=q^{\beta}\gap^{-1}$).

Then we can write
%
%
\begin{eqnarray}\label
{eq:deuxtermssomme}
\mathbb{E}_{}\bigl[{X_t^2}\bigr]&=&
\mathbb{E}_{}\Biggl[{ \Biggl(\sum_{k=1}^{t/\tau
}X_{k\tau}-X_{(k-1)\tau }
\Biggr)^2}\Biggr]
\nonumber
\\
&=&\sum_{k=1}^{t/\tau}\mathbb{E}_{}
\bigl[{(X_{k\tau}-X_{(k-1)\tau} )^2}\bigr]
\nonumber
\\[-8pt]
\\[-8pt]
\nonumber
&&{}+\sum
_{k\neq k'}\mathbb{E}_{}\bigl[{(X_{k\tau}-X_{(k-1)\tau}
) (X_{k'\tau
}-X_{(k'-1)\tau} )}\bigr]
\\
&=&\frac{t}{\tau}\mathbb{E}_{}\bigl[{X_{\tau}^2}
\bigr]+\sum_{k\neq k'}\mathbb{E}_{}
\bigl[{(X_{k\tau }-X_{(k-1)\tau} ) (X_{k'\tau
}-X_{(k'-1)\tau}
)}\bigr]\nonumber.
\end{eqnarray}
We need to show that \eqref{eq:deuxtermssomme} is smaller than
$tq^{-\alpha}\gap$ for some $\alpha$ when $\tau$ is well chosen. We are
going to bound the first term using the fact that energy barriers make
it very costly to cross a distance greater than $1/q$ in time $\tau
\lesssim\gap^{-1}$. To bound the second term, we use the symmetry of
the model and the fact that the process seen from the tracer has a
positive spectral gap.

%
\begin{prop}\label{prop:Xtau2}
There exists $\beta,C<\infty$ such that, if $\tau= q^{\beta}\gap^{-1}$,
%
%
\begin{equation}
\mathbb{E}_{}\bigl[{X_{\tau}^2}\bigr]\leq
Cq^{-C}.
\end{equation}
\end{prop}

First we need two lemmas that rely on precise estimates on the spectral
gap of the East model on lengths of order at most $1/q$, and related
energy barriers, that have been established in \citet
{timescaleseparation}. We start by showing a precise comparison between
the relaxation time in infinite volume and the relaxation time in
volume $1/q$. Recall that it was shown in \citet{KCSM} that for any
$\delta>0$,
%
%
\begin{equation}
\label{eq:gapeastdelta} \gap^{-1}\leq C_\delta \biggl(
\frac{1}{q} \biggr)^{\log_2(1/q)/(2-\delta)}.
\end{equation}

%
\begin{lemme}\label{lem:newgapeast}
Let $n= \lceil\log_2(1/q) \rceil$ and $T_{\mathrm{rel}}(L)$ be the
relaxation time of the East model on length $L$ with empty boundary
condition. Then there exist finite constants $C,C'$ such that
%
%
\begin{equation}
\gap^{-1}\leq Cq^{-C}T_{\mathrm{rel}}(1/q)\leq
C'q^{-C'}\frac{n!}{q^n2^{\bigl({n\atop2}\bigr)}}.
\end{equation}
\end{lemme}

\begin{pf}
The second inequality follows immediately from Theorem~2 in \citet
{timescaleseparation}. To prove the first one, we refine the bisection
technique used in \citet{KCSM} to prove \eqref{eq:gapeastdelta}. Let
$\delta(q)=10/\log(1/q)$, $l_k=2^k$, $\delta_k= \lfloor
l_k^{1-\delta/2} \rfloor$, $s_k= \lfloor l_k^{\delta/6}
\rfloor$. These are the same definitions as in \citet{KCSM}, except that
instead of a fixed $\delta>0$, we take $\delta$ to $0$ with~$q$. With
these definitions, we have for every $k\geq k_\delta:=6/\delta$ the
following estimate\footnote{This condition is not necessary, but
sufficient; it comes from the fact that Lemma~4.2 in \citet{KCSM} has to
be satisfied in order to apply the bisection technique.} [see (6.3)
in~\citet{KCSM}]:
%
%
\begin{equation}
\label{eq:bisection} \gap^{-1}\leq T_{\mathrm{rel}} \bigl(l_k+l_k^{1-\delta/6}
\bigr)\prod_{j=k}^\infty \biggl(
\frac{1}{1-p^{\delta_j/2}} \biggr)\prod_{j=k}^\infty
\bigl(1+s_j^{-1} \bigr).
\end{equation}
As in \citet{KCSM}, let
%
%
\begin{equation}
\label{eq:j*} j_*=\min \bigl\{j | p^{\delta_j/2}\leq e^{-1} \bigr\}
\approx\log_2(1/q)/(1-\delta/2).
\end{equation}

As long as $j_*\geq k_\delta$, which is true thanks to our choice of
$\delta$, we can replace $k$ by $j_*$ in \eqref{eq:bisection}. Now we
have [see the computations in \citet{KCSM}, top of page 484 for the
first estimate]
%
%
\begin{eqnarray}
\prod_{j=j_*}^\infty \biggl(\frac{1}{1-p^{\delta_j/2}}
\biggr)&\leq& C,
\\
\prod_{j=j_*}^\infty \bigl(1+s_j^{-1}
\bigr) &\leq& q^{-C},
\end{eqnarray}
for $C$ some constant not depending on $q$. Noticing that
$l_{j_*}+l_{j_*}^{1-\delta/6}\leq d/q$ for some constant $d$, we get
%
%
\begin{equation}
\label{eq:newgapestimate} \gap^{-1}\leq Cq^{-C}T_{\mathrm{rel}}(d/q).
\end{equation}
Now it is enough to recall Theorem~4 in \citet{timescaleseparation},
that states that there is no time scale separation on scale $1/q$
%
%
\begin{eqnarray}
T_{\mathrm{rel}}(d/q)&\sim&T_{\mathrm{rel}}(1/q)\label
{eq:equivalenceoftimes}.
\end{eqnarray}
\upqed\end{pf}

Now we can use Lemma~\ref{lem:newgapeast} to prove the following
estimate, which basically means that in times smaller than $\gap^{-1}$,
it will be extremely difficult for the system to erase a row of $1/q$ ones.

%
%
\begin{lemme}\label{lem:energybarriers}
Recall that $\tau= q^{\beta}\gap^{-1}$. Let $l= 1/q$ and $\mathbb
{P}_{\mathbf {1}0}({\cdot})$ denote (abusively) the law of the East
process on
$\mathbb
{Z}$ starting from a configuration equal to one on $\lbrace
1,\ldots,l\rbrace$, with a zero in $l+1$. Let $T_0$ be the first time
there is a zero at $1$. Independently of the choice of the initial
configuration outside $\lbrace1,\ldots,l,l+1\rbrace$, we have, if
$\beta$
is large enough (independently of $q$),
%
%
\begin{equation}
\label{eq:energybarrier} \mathbb{P}_{\mathbf{1}0}({T_0\leq\tau })\leq Cq.
\end{equation}
\end{lemme}

\begin{pf}
In \citet{timescaleseparation},\hskip.2pt\footnote{Note that the orientation
convention is reversed in this paper: contrary to this paper, the
constraint that has to be satisfied to update $x$ is that $x-1$ should
be empty.} the authors define a certain set $\partial A_*$ of
configurations in $\lbrace0,1\rbrace^l$ that has two interesting
properties [it is defined in paragraph 5.2.1 of \citet
{timescaleseparation}, the properties below are stated in Remark~5.8
and Corollary~5.10]:
\begin{itemize}
\item Starting from a configuration equal to one on $\lbrace
1,\ldots,l\rbrace$, with a zero in $0$, in order to put a one in $0$
before time $\tau$, the dynamics restricted to $\lbrace1,\ldots,l\rbrace$
has to go through the set $\partial A_*$ at some time $s\leq\tau$.
\item For some $\alpha'<\infty$, if $n=\lceil\log_2 l\rceil$,
%
%
\begin{equation}
\label{eq:poidsA*} \mu(\partial A_*)\leq\frac{q^n2^{\bigl({n\atop
2}\bigr)}}{n!}q^{-\alpha'}.
\end{equation}
\end{itemize}
Put another way, $\partial A_*$ is a bottleneck separating the events
$ \lbrace\eta_0=\eta_{l+1}=0, \eta_1=\cdots=\eta_l=1 \rbrace$
and \mbox{$\lbrace\eta_0=1\rbrace$} in the East dynamics.

Call $\tau_0$ the first time there is a one in $0$. Denote (abusively)
by $0\mathbf{1}0$ any configuration equal to zero in $0$ and $l+1$, and
to one on $\lbrace1,\ldots,l\rbrace$, by $T$ an exponential variable of
parameter $2$ independent of $T_0$, and by $\tau_0$ the first time at
which there is a one in position $0$. Notice that once there is a zero
in $1$, if the clock attached to site $0$ rings before that attached to
$1$, and if the associated Bernoulli variable is a one, then the
configuration at site $0$ takes value one. So that
%
%
\begin{equation}
\label{eq:T_0+T} \frac{1-q}{2}\mathbb{P}_{\mathbf{1}0}({T_0+T
\leq \tau+1/2}) \leq\mathbb{P}_{0\mathbf{1}0}({\tau_0\leq\tau+1/2}),
\end{equation}
where $\mathbf{1}0$ and $0\mathbf{1}0$ are equal except maybe in $0$.
The constant $1/2$ appears to allow the following estimate:
%
%
\begin{eqnarray}
\label{eq:T_0tau_0} \mathbb{P}_{\mathbf{1}0}({T_0+T\leq\tau +1/2})&\geq&
\mathbb{P}_{\mathbf{1}0}({T_0\leq\tau})\mathbb {P}_{}({T
\leq1/2})
\nonumber
\\[-8pt]
\\[-8pt]
\nonumber
&=& \bigl(1-e^{-1} \bigr)\mathbb{P}_{\mathbf{1}0}({T_0
\leq\tau}).
\end{eqnarray}
Equations \eqref{eq:T_0+T} and \eqref{eq:T_0tau_0} yield
%
%
\begin{equation}
\label{eq:T0leqtau} \mathbb{P}_{\mathbf{1}0}({T_0\leq\tau})\leq
\frac
{2}{(1-q)(1-e^{-1})} \mathbb{P}_{0\mathbf{1}0}({\tau_0\leq\tau+1/2}).
\end{equation}
Now we use the first property of $\partial A_*$ to get
%
%
\begin{eqnarray}
\mathbb{P}_{0\mathbf{1}0}({\tau_0\leq\tau+1/2})&\leq&\mathbb
{P}_{\mathbf{1}0}\bigl({\exists s\leq\tau+1/2 \mbox{ s.t. } \bigl(\omega(s)
\bigr)_{[1,l]}\in \partial A_*}\bigr)\label{eq:passageparbottleneck}.
\end{eqnarray}
To evaluate the RHS, we condition on $N_{\tau+1/2}$ the number of rings
occurring in $[1,l]$ before time $\tau+1/2$ in the graphical
construction with a union bound to get
%
%
\begin{eqnarray}
\label{eq:taumu(A*)}&& \mathbb{P}_{\mathbf{1}0}\bigl({\exists s\leq\tau +1/2 \mbox
{ s.t. } \bigl(\omega (s) \bigr)_{[1,l]}\in\partial A_*}\bigr)\nonumber\\
&&\qquad\leq
\mathbb{E}_{}[{N_{\tau+1/2}}]\sup_{s\leq
\tau+1/2}
\mathbb{P}_{\mathbf{1}0}\bigl({ \bigl(\omega(s) \bigr)_{[1,l]} \in
\partial A_*}\bigr)
\nonumber
\\
&&\qquad\leq (\tau+1/2 )l\sum_{\sigma\in\partial A_*}\sup
_{s\leq
\tau+1/2}\mathbb{P}_{\mathbf{1}0}\bigl({ \bigl(\omega(s)
\bigr)_{[1,l]}= \sigma}\bigr)
\nonumber
\\[-8pt]
\\[-8pt]
\nonumber
&&\qquad\leq (\tau+1/2 )l\sum_{\sigma\in\partial A_*}p^{-l}\mu(
\sigma)
\\
&&\qquad\leq (\tau+1/2 )l(1-q)^{-l}\mu(\partial A_*)
\nonumber
\\
&&\qquad\leq (\tau+1/2 ) (1-q)^{-l}\frac{q^n2^{\bigl({n\atop
2}\bigr)}}{n!}q^{-(\alpha'+1)},\nonumber
\end{eqnarray}
where we used \eqref{eq:distinguishedestimate} with the distinguished
zero starting at $l+1$ to get the third inequality, and the second
property of $\partial A_*$ \eqref{eq:poidsA*} to get the last one. Now
collect~\eqref{eq:T0leqtau}, \eqref{eq:passageparbottleneck} and
\eqref
{eq:taumu(A*)} to get
%
%
\begin{equation}
\mathbb{P}_{\mathbf{1}0}({T_0\leq\tau})\leq\frac
{2}{(1-q)(1-e^{-1})} (
\tau +1/2 ) (1-q)^{-l}\frac{q^n2^{\bigl({n\atop2}\bigr)}}{n!}q^{-(\alpha'+1)}.
\end{equation}
For $q$ small enough and $\tau= q^{\beta}\gap^{-1}$, $\tau+1/2\leq
q^{\beta-1}\gap^{-1}$, so that\break Lemma~\ref{lem:newgapeast} yields
%
%
\begin{equation}
\mathbb{P}_{\mathbf{1}0}({T_0\leq\tau})\leq Cq^{-\alpha''}q^{\beta-1},
\end{equation}
for some $C,\alpha''$ independent of $q$.
\end{pf}

%
\begin{rem}
An anonymous referee suggested an alternative proof for this lemma,
relying directly on Proposition~3.2 and Theorem~1 in \citet
{timescaleseparation} and Lemma~\ref{lem:newgapeast}, outlined as
follows. $\mathbb{P}_{0\mathbf{1}0}({\tau_0\leq t})\leq
et/T_{\mathrm{hit}}(1/q)$ by (3.3)
in \citet{timescaleseparation}, and by Theorem~1 in \citet
{timescaleseparation}, $T_{\mathrm{hit}}(1/q)\geq c T_{\mathrm{rel}}(1/q)$. Lemma~\ref
{lem:newgapeast} then yields the conclusion. In order to carry this
(more efficient) proof rigorously, one would just need to check that
the above results can be extended to infinite volume dynamics with
distinguished zero starting at $l+1$ (or guarantee that the initial
zero at $l+1$ has not moved by time $\tau$). We keep the above proof in
order to evidence the role of the energy barriers involved in confining
the tracer, the most relevant part of the proof in that respect being
the set of equations~(\ref{eq:taumu(A*)}).
\end{rem}

\begin{pf*}{Proof of Proposition~\ref{prop:Xtau2}}
First of all, let us reformulate what we want to show.
%
%
\begin{eqnarray}\label{eq:Xtaucarre}
\mathbb{E}_{}\bigl[{X_{\tau}^2}\bigr]&=&\sum
_{x=1}^\infty(2x-1) \mathbb{P}_{}\bigl({|X_\tau|
\geq x}\bigr)
\nonumber
\\
&=&2\sum_{x=1}^\infty(2x-1)
\mathbb{P}_{}({X_\tau\geq x})
\\
&\leq&4\sum_{m=1}^\infty q^{-m}
\mathbb{P}_{}\bigl({X_\tau\geq q^{-m}}
\bigr).\nonumber
\end{eqnarray}
In light of Lemma~\ref{lem:energybarriers}, we can now notice that in
order to have $X_\tau\geq q^{-m}$ for $m\geq2$, the system will have
to overcome a large number of energy barriers (i.e., rows of ones of
length larger than $1/q$), so that the probability of this event will
become very small.

Fix $m>2$, and let us study $\mathbb{P}_{}({X_\tau\geq q^{-m}})$.
Throughout the
proof, to simplify the notation, if $C(q)$ is a quantity going to
infinity when $q\rightarrow0$, we will not make the distinction
between $C(q)$ and $\lfloor C(q)\rfloor$. We divide $\lbrace
0,\ldots,q^{-m}\rbrace$ into $q^{-m+2}(3m)^{-1}$ groups of $3m $
blocks of
length $q^{-2}$. Given a configuration, we say that a block of $q^{-2}$
sites is well behaved if we can find a row of consecutive ones of
length at least $1/q$ that ends with a zero inside it. We can estimate
the probability of a block having this property by
%
%
\begin{equation}
\label{eq:probagentil} \mu(\mbox{a given block is not well-behaved})\leq
\bigl(1-q(1-q)^{1/q} \bigr)^{1/q}\leq c<1
\end{equation}
for some constant $c$.

Let $A$ be the event that in all of these $q^{-m+2}(3m)^{-1}$ groups of
blocks, there is one of the $3m$ blocks that is not well behaved. With
this definition, on $A^c$, there is a group of $3m$ well behaved
blocks. Let us estimate the probability of $A$ under $\mu$ using
\eqref
{eq:probagentil}
%
%
\begin{eqnarray}
\label{eq:probaA2} \mu(A)&\leq& \bigl(1-\mu(\mbox{a given block is
well-behaved})^{3m} \bigr)^{q^{-m+2}(3m)^{-1}}
\nonumber
\\
&\leq& \bigl(1-(1-c)^{3m} \bigr)^{q^{-m+2}(3m)^{-1}}
\\
&\leq&e^{-Cq^2q^{-\gamma m}},\nonumber
\end{eqnarray}
with $\gamma>0$, $C<\infty$.

So we can write
%
%
\begin{equation}
\label{eq:probatracereast} \mathbb{P}_{}\bigl({X_\tau\geq
q^{-m}}\bigr)\leq \mu(A)+\mu \bigl( \mathbf{1}_{A^c}(\eta)
\mathbb{P}_{\eta}\bigl({X_\tau\geq q^{-m}}\bigr)
\bigr).
\end{equation}
Denote by $B_1$ the first block of length $q^{-2}$,
$B_2=B_1+q^{-2},\ldots, B_{3m}=B_1+ (3m-1 )q^{-2}$. We have the
following estimate:
%
%
\begin{eqnarray}
\label{eq:Xtau3} &&\mu \bigl(\mathbf{1}_{A^c}(\eta)\mathbb{P}_{\eta
}
\bigl({X_\tau\geq q^{-m}}\bigr) \bigr)
\nonumber
\\[-8pt]
\\[-8pt]
\nonumber
&&\qquad\leq
q^{-m+2}(3m)^{-1} \mu \Biggl(\prod
_{i=1}^{3m} \mathbf{1}_{B_i\ \mathrm{well\mbox{-}behaved\ }}(\eta)
\mathbb{P}_{\eta}\bigl({X_\tau\geq q^{-m}}\bigr)
\Biggr).
\end{eqnarray}
Let $\eta$ be a configuration in which all the $B_i$ are well behaved.
Let $x_i$ be the starting point of the first row of $1/q$ ones ended by
a zero in $B_i$, and $T_i$ the first time this site is empty. We denote
by $ (\xi_i(s) )_{s\leq\tau}$ the trajectory of the
distinguished zero started from the position of the zero at the end of
the row of ones starting at $x_i$, up to time~$\tau$:
%
%
\begin{eqnarray}
\label{eq:Xtau1} &&\mathbb{P}_{\eta}\bigl({X_\tau\geq
q^{-m}}\bigr)\nonumber\\
&&\qquad\leq \mathbb{P}_{\eta}({\forall i=1,\ldots,3m\
T_i\leq\tau})
\nonumber
\\
&&\qquad\leq \mathbb{P}_{\eta}({T_{3m}\leq\tau})
\mathbb{P}_{\eta
}({\forall i=1,\ldots,3m-1\ T_i\leq\tau|
T_{3m}\leq\tau})
\\
&&\qquad\leq \mathbb{P}_{\eta}({T_{3m}\leq\tau})
\nonumber\\
&&\qquad\quad{}\times\mathbb{E}_{\eta
}\bigl[{\mathbb{P}_{\eta}\bigl({\forall i=1,
\ldots,3m-1\ T_i\leq\tau| \bigl(\xi_{3m-1}(s)
\bigr)_{s\leq\tau}}\bigr) | T_{3m}\leq\tau}\bigr]\nonumber,
\end{eqnarray}
since the dynamics on the left of $x_{3m-1}+1/q$ knowing $ (\xi
_{3m-1}(s) )_{s\leq\tau}$ does not depend on what happens on the
right of $ (\xi_{3m-1}(s) )_{s\leq\tau}$.

Let us show iteratively that, uniformly in the trajectory $ (\xi
_{k}(s) )_{s\leq\tau}$,
%
%
\begin{equation}
\label{eq:Xtau2} \mathbb{P}_{\eta}\bigl({\forall i=1,\ldots, k\
T_i\leq \tau| \bigl( \xi_{k}(s) \bigr)_{s\leq\tau}}
\bigr)\leq(Cq )^{k}.
\end{equation}
For $k=1$, mutatis mutandis, the proof of Lemma~\ref
{lem:energybarriers} applies. Let $k>1$.

$\mathbb{P}_{\eta}({\forall i=1,\ldots, k\ T_i\leq\tau| (\xi
_{k}(s) )_{s\leq\tau}})$ is also
%
%
\begin{equation}\qquad
\mathbb{P}_{\eta}\bigl({T_k\leq\tau| \bigl(
\xi_{k}(s) \bigr)_{s\leq
\tau}}\bigr)\mathbb{P}_{\eta}
\bigl({\forall i=1,\ldots, k-1\ T_i\leq\tau | \bigl(
\xi_{k}(s) \bigr)_{s\leq \tau},T_k\leq\tau}\bigr),
\end{equation}
which can be rewritten
%
%
\begin{eqnarray}\qquad
&&\mathbb{P}_{\eta}\bigl({T_k\leq\tau| \bigl(
\xi_{k}(s) \bigr)_{s\leq
\tau}}\bigr)
\nonumber
\\[-8pt]
\\[-8pt]
\nonumber
&&\qquad{}\times\mathbb{E}_{\eta}
\bigl[{\mathbb{P}_{\eta}\bigl({\forall i=1,\ldots, k-1\ T_i
\leq\tau| \bigl(\xi_{k-1}(s) \bigr)_{s\leq\tau}}\bigr) | \bigl(
\xi_{k}(s) \bigr)_{s\leq\tau },T_k\leq\tau}\bigr],
\end{eqnarray}
and the induction hypothesis applies.

Putting together \eqref{eq:Xtau1}, \eqref{eq:Xtau2} and \eqref
{eq:Xtau3}, we get for some constant $C$
%
%
\begin{equation}
\mu \bigl(\mathbf{1}_{A^c}(\eta)\mathbb{P}_{\eta}
\bigl({X_\tau\geq q^{-m}}\bigr) \bigr)\leq(Cq)^{2m}.
\end{equation}
Recalling \eqref{eq:probatracereast}, \eqref{eq:probaA2} and \eqref
{eq:Xtaucarre}, we get Proposition~\ref{prop:Xtau2}.
\end{pf*}

What now remains is to show there is enough decorrelation to bound the
second sum in \eqref{eq:deuxtermssomme}. This is not difficult, once we
make the following remark.

%
\begin{lemme}\label{lem:spectralgaps}
Denote by $\gap_T$ the spectral gap of the process seen from the tracer
[recall \eqref{eq:processseenfromtracer}]
%
%
\begin{equation}
\gap_T=\inf\frac{-\mu(f\mathcal{L}f)}{\operatorname{Var}_\mu(f)},
\end{equation}
where the infimum is taken over nonconstant functions $f\in L^2(\mu)$.
Then we have
%
%
\begin{equation}
\label{eq:spectralgaps} \gap_T\geq\gap.
\end{equation}
\end{lemme}

\begin{pf}
This follows directly from \eqref{eq:var2} and the definition of $\gap$
and $\gap_T$ [recall \eqref{eq:defgap}].
\end{pf}

Now we are armed to study the terms $\mathbb{E}_{}[{(X_{k\tau
}-X_{(k-1)\tau } ) (X_{k'\tau}-X_{(k'-1)\tau} )}]$. First of all, by
stationarity, this quantity depends only on $\tau$ and $|k-k'|$.
Therefore, we only need to study $\mathbb{E}_{}[{X_{\tau} (X_{k\tau
}-X_{(k-1)\tau} )}]$ for $k\geq2$. In fact, using the
Cauchy--Schwarz inequality and Proposition~\ref{prop:Xtau2}, we only
need to study this term for $k\geq3$, which allows some decorrelation
to take place between times $\tau$ and $(k-1)\tau$. Let us denote by
$ (P^T_s )_{s\geq0}$ the semigroup associated to $\mathcal
{L}$. $\mathbb{E}_{(\omega,x)}[{\cdot}]$ will denote the law of the
process with
generator $\mathcal{L}_0$ starting from the configuration $\omega$ with
the tracer in position $x$ ($\mathbb{E}_{}[{\cdot}]$ is still the law
of the
process starting from $\mu$ and the tracer at the origin). Using
successively the Markov property at time $\tau$, we can write
%
%
\begin{eqnarray}
\mathbb{E}_{}\bigl[{X_{\tau} (X_{k\tau}-X_{(k-1)\tau}
)}\bigr]&=&\mathbb {E}_{}\bigl[{X_\tau
\mathbb{E}_{ (\omega(\tau),X_\tau
)}\bigl[{X'_{(k-1)\tau}-X'_{(k-2)\tau}}
\bigr]}\bigr],
\end{eqnarray}
where $ (X'_s )_{s\geq0}$ denotes the trajectory of the
tracer under the law $\mathbb{E}_{(\omega(\tau),X_\tau)}[{\cdot
}]$. Now we use
successively the Cauchy--Schwarz inequality and stationarity of the
process seen from the tracer to get
%
%
\begin{eqnarray}
\mathbb{E}_{}\bigl[{X_{\tau} (X_{k\tau}-X_{(k-1)\tau}
)}\bigr]^2& \leq&\mathbb{E}_{}\bigl[{X_\tau^2}
\bigr]\mathbb{E}_{}\bigl[{\mathbb{E}_{
(\omega( \tau),X_\tau )}[{X_{(k-1)\tau}-X_{(k-2)\tau }}]^2}
\bigr]
\nonumber
\\
&\leq&\mathbb{E}_{}\bigl[{X_\tau^2}\bigr]
\mathbb{E}_{}\bigl[{\mathbb{E}_{ (
(\omega(\tau)  )_{X_\tau +\cdot},0
)}[{X_{(k-1)\tau}-X_{(k-2)\tau}}]^2}
\bigr]
\\
&\leq&\mathbb{E}_{}\bigl[{X_\tau^2}\bigr]\mu
\bigl(\mathbb{E}_{(\omega,0
)}[{X_{(k-1)\tau }-X_{(k-2)\tau}}]^2
\bigr).\nonumber
\end{eqnarray}
Let us focus on $\mathbb{E}_{(\omega,0 )}[{X_{(k-1)\tau
}-X_{(k-2)\tau }}]$. Using the Markov property at time $(k-2)\tau$,
we get
%
%
\begin{eqnarray}
\mathbb{E}_{(\omega,0 )}[{X_{(k-1)\tau}-X_{(k-2)\tau}}]&=&\mathbb
{E}_{(\omega,0 )}\bigl[{\mathbb{E}_{ (\omega ((k-2)\tau
),X_{(k-2)\tau}  )}\bigl[{X'_\tau
-X'_0}\bigr]}\bigr]
\nonumber
\\
&=&\mathbb{E}_{(\omega,0 )}\bigl[{\mathbb{E}_{ (  (\omega
((k-2)\tau  )  )_{X_{(k-2)\tau+\cdot},0}  )}
\bigl[{X'_\tau}\bigr]}\bigr]
\\
&=&P_{(k-2)\tau}^Tg(\omega),\nonumber
\end{eqnarray}
where $g(\omega)=\mathbb{E}_{(\omega,0)}[{X_\tau}]$, and the $X'_s$
in the first
and second line denote respectively the trajectory of the tracer under
the laws $\mathbb{E}_{(\omega((k-2)\tau),X_{(k-2)\tau} )}[{\cdot
}]$ and
$\mathbb{E}_{( (\omega((k-2)\tau) )_{X_{(k-2)\tau+\cdot}},0
)}[{\cdot}]$. Therefore, using the spectral gap inequality and the fact
that $g$ is a mean-zero function in $L^2(\mu)$ thanks to stationarity
and Proposition~\ref{prop:Xtau2}, we get
%
%
\begin{eqnarray}
\mathbb{E}_{}\bigl[{(X_{\tau} ) (X_{k\tau}-X_{(k-1)\tau}
)}\bigr]^2&\leq &\mathbb{E}_{}\bigl[{X_\tau^2}
\bigr]\mu \bigl( \bigl(P_{(k-2)\tau}^Tg \bigr)^2
\bigr)
\nonumber
\\
&\leq& \mathbb{E}_{}\bigl[{X_\tau^2}
\bigr]^2e^{-2(k-2)\tau\gap_T}
\\
&\leq& \mathbb{E}_{}\bigl[{X_\tau^2}
\bigr]^2e^{-2(k-2)q^{\beta}}.\nonumber
\end{eqnarray}

Since $\sum_{k\geq1}e^{-kq^\beta}\lesssim q^{-\beta}$, the second term
in \eqref{eq:deuxtermssomme} is
%
%
\begin{equation}
\sum_{k\neq k'}\mathbb{E}_{}
\bigl[{(X_{k\tau}-X_{(k-1)\tau} ) (X_{k'\tau }-X_{(k'-1)\tau}
)}\bigr]\leq C \lfloor t/\tau\rfloor\mathbb {E}_{}
\bigl[{X_\tau^2}\bigr]q^{-\beta}.
\end{equation}
Putting this into \eqref{eq:deuxtermssomme} together with
Proposition~\ref{prop:Xtau2}, we get Theorem~\ref{th:east}.
\end{pf*}

\begin{appendix}\label{app}
\section*{Appendix: An alternative proof in the FA-1f model}

When the environment is given by the one-spin Fredrickson--Andersen
model (FA-1f), in which $c_x(\eta)=1-\prod_{i=1}^d\eta_{e_i}\eta
_{-e_i}$ (the constraint requires at least one nearest neighbor to be
empty), the diffusion coefficient at low density is of order~$q^2$.
This means that in this particular case, the correct order is already
given by the first term in~\eqref{eq:formuleD}, which allows us to
design another strategy to find the lower bound in Theorem~\ref
{th:kzeros} when $k=1$. Since the diffusion coefficient is of order
lower than $q^2$ in the $k$-zeros model with $k> 1$, this technique
does not apply. For simplicity, we write the proof in dimension $d=1$.

We follow the strategy devised to prove Lemma~6.25 in \citet
{ollalandim}; that is, we prove that
%
%
\begin{equation}
\label{eq:alternativeproof} \sup \bigl\{2\mu(jf)-\mathcal{D}(f) \bigr\}\leq cq^2,
\end{equation}
where $c<1$ does not depend on $q$ and $\mathcal{D}(f)=-\mu
(f\mathcal{L}f )$. Seeing \eqref{eq:formuleD} and \eqref
{eq:varformula1}, this is sufficient to prove Theorem~\ref{th:kzeros}
when $k=1$, $d=1$. To obtain that result, we define [recall \eqref{eq:defrate}]
%
%
\begin{eqnarray}
\mathcal{D}_{\mathrm{jump}}(f)&=&\frac{1}{2} \mu \bigl(\bar{
\eta}_0\bar{\eta}_{\alpha e_i} \bigl[f(\eta_{\alpha e_i+\cdot})-f(
\eta) \bigr]^2 \bigr),
\\
\mathcal{D}_{\mathrm{FA}}(f)&=&\frac{1}{2}\sum
_{y\in\Z}\mu \bigl(c_y(\eta)r_y(\eta)
\bigl[f \bigl(\eta^y \bigr)-f(\eta) \bigr]^2 \bigr),
\end{eqnarray}
so that $\mathcal{D}(f)=\mathcal{D}_{\mathrm{jump}}(f)+\mathcal{D}_{\mathrm{FA}}(f)$, and
we show separately that for all $f$,
%
%
\begin{eqnarray}
2\mu(jf)-\mathcal{D}_{\mathrm{jump}}(f)&\leq&q^2,\label{eq:bornejump}
\\
2\mu(jf)-\mathcal{D}_{\mathrm{FA}}(f)&\leq&Cq^2,\label{eq:borneFA}
\end{eqnarray}
where $C\geq1$ is a constant that does not depend on $q$. To get the
result from \eqref{eq:bornejump} and \eqref{eq:borneFA}, we write that
for any $\lambda>0$, for any local function $f$.
\begin{eqnarray*}
&&\lambda^{-1} \bigl(2\mu(jf)-\mathcal{D}_{\mathrm{jump}}(f)-\mathcal
{D}_{\mathrm{FA}}(f) \bigr)\\
&&\qquad=2\mu \bigl(j\lambda^{-1}f \bigr)-\lambda
\mathcal{D}_{\mathrm{jump}} \bigl(\lambda^{-1}f \bigr)-\lambda
\mathcal{D}_{\mathrm{FA}} \bigl(\lambda^{-1}f \bigr)
\end{eqnarray*}
so that
\[
\lambda^{-1}\sup \bigl\{2\mu(jf)-\mathcal{D}_{\mathrm{jump}}(f)-
\mathcal{D}_{\mathrm{FA}}(f) \bigr\}\leq\sup \bigl\{2\mu(jg )-\lambda
\mathcal{D}_{\mathrm{jump}} (g )-\lambda\mathcal{D}_{\mathrm{FA}} (g ) \bigr\}.
\]
Take, for instance, $\lambda=C/(C+1)$. We have $\lambda\geq1-\lambda$,
so that
\begin{eqnarray*}
&&\lambda^{-1}\sup \bigl\{2\mu(jf)-\mathcal{D}_{\mathrm{jump}}(f)-
\mathcal{D}_{\mathrm{FA}}(f) \bigr\}\\
&&\qquad\leq \sup \bigl\{2\mu(jg )-\lambda
\mathcal{D}_{\mathrm{jump}} (g )-(1-\lambda)\mathcal{D}_{\mathrm{FA}} (g ) \bigr
\}
\\
&&\qquad\leq \bigl[\lambda+(1-\lambda)C \bigr]q^2=q^2,
\end{eqnarray*}
using \eqref{eq:bornejump} and \eqref{eq:borneFA}, so that \eqref
{eq:alternativeproof} is proven.

(1) Proof of \eqref{eq:bornejump}.

For any local function $f$, we can rewrite $\mu(jf)$ in terms of the
``jumps'' $\eta\rightarrow\eta_{1+\cdot}$ and $\eta\rightarrow
\eta_{-1+}$
\[
2\mu(jf)=-\mu \bigl(\bar{\eta}_0\bar{\eta}_1
\bigl[f(\eta_{1+\cdot
})-f(\eta) \bigr] \bigr)+\mu \bigl(\bar{
\eta}_0\bar{\eta}_{-1} \bigl[f(\eta_{-1+\cdot})-f(
\eta) \bigr] \bigr).
\]
Now using the inequality $ab\leq(a^2+b^2)/2$, the Dirichlet form
$\mathcal{D}_{\mathrm{jump}}(f)$ appears on the RHS,
\begin{eqnarray*}
2\mu(jf)&\leq&q^2+\tfrac{1}{2}\mu \bigl(\bar{
\eta}_0\bar{\eta}_1 \bigl[f(\eta_{1+\cdot})-f(
\eta) \bigr]^2 \bigr)+\tfrac{1}{2}\mu \bigl(\bar{\eta
}_0\bar{\eta}_{-1} \bigl[f(\eta_{-1+\cdot})-f(\eta)
\bigr]^2 \bigr)
\\
&\leq&q^2+\mathcal{D}_{\mathrm{jump}}(f).
\end{eqnarray*}

(2) Proof of \eqref{eq:borneFA}.

We need only to prove it for small $q$. First we make a few
computations to express $\mu(jf)$ in terms of allowed flips ($\eta
\rightarrow\eta^1$ or $\eta\rightarrow\eta^{-1}$). Then we use the same
optimization technique performed in the proof of Lemma~6.13 in \citet
{ollalandim} to get the desired bound. We have the following equalities:
%
%
\begin{eqnarray}
\mu \bigl(\bar{\eta}_0\bar{\eta}_1f(\eta) \bigr)&=&
\frac{q}{1-2q}\mu \bigl(\bar{\eta}_0 \bigl[f \bigl(
\eta^1 \bigr)-f(\eta) \bigr] \bigr)+\frac{q}{1-q}\mu \bigl(\bar{
\eta}_0\eta_1f(\eta) \bigr),
\\
\qquad\quad\mu \bigl(\bar{\eta}_0\bar{\eta}_{-1}f(\eta) \bigr)&=&
\frac{q}{1-2q}\mu \bigl(\bar{\eta}_0 \bigl[f \bigl(
\eta^{-1} \bigr)-f(\eta) \bigr] \bigr)+\frac
{q}{1-q}\mu \bigl(\bar{
\eta}_0\eta_{-1}f(\eta) \bigr),
\\
\mu \bigl(\bar{\eta}_0\eta_1f(\eta) \bigr)&=&(1-q)\mu
\bigl(\bar{\eta}_0\bar{\eta}_1 \bigl[f \bigl(
\eta^1 \bigr)-f(\eta) \bigr] \bigr)+(1-q)\mu \bigl(\bar{
\eta}_0f(\eta) \bigr),
\\
\mu \bigl(\bar{\eta}_0\eta_{-1}f(\eta) \bigr)&=&(1-q)\mu
\bigl(\bar{\eta}_0\bar{\eta}_{-1} \bigl[f \bigl(
\eta^{-1} \bigr)-f(\eta) \bigr] \bigr)+(1-q)\mu \bigl(\bar{
\eta}_0f(\eta) \bigr).
\end{eqnarray}
So that, computing differences, we get
%
%
\begin{eqnarray}
\mu(jf)&=&\frac{q}{p-q} \bigl[\mu \bigl(\bar{\eta}_0 \bigl[f
\bigl(\eta^1 \bigr)-f(\eta) \bigr] \bigr)-\mu \bigl(\bar{
\eta}_0 \bigl[f \bigl(\eta^{-1} \bigr)-f(\eta) \bigr]
\bigr) \bigr]
\nonumber
\\[-8pt]
\\[-8pt]
\nonumber
&&{}+ q \bigl[\mu \bigl(\bar{\eta}_0\bar{\eta}_1 \bigl[f
\bigl(\eta^1 \bigr)-f(\eta) \bigr] \bigr)-\mu \bigl(\bar{
\eta}_0\bar{\eta}_{-1} \bigl[f \bigl(\eta^{-1}
\bigr)-f(\eta) \bigr] \bigr) \bigr].
\nonumber
\end{eqnarray}
Assume $q<1/2$. Using the inequality $ab\leq(a^2+b^2)/2$, we get for
any $\alpha,\beta>0$,
\begin{eqnarray*}
&&\frac{\mu(jf)}{q}\\
&&\qquad\leq \frac{1}{1-2q} \biggl\{\alpha q+\frac
{1}{2\alpha
}
\bigl[\mu \bigl(\bar{\eta}_0 \bigl[f \bigl(\eta^1
\bigr)-f(\eta) \bigr]^2 \bigr)+\mu \bigl(\bar{\eta}_0
\bigl[f \bigl(\eta^{-1} \bigr)-f(\eta) \bigr]^2 \bigr)
\bigr] \biggr\}\\
&&\qquad\quad{}+\beta
\nonumber
\\
&&\qquad\quad{}+ \frac{1}{2\beta} \bigl[\mu \bigl(\bar{\eta}_0\bar{
\eta}_1 \bigl[f \bigl(\eta^1 \bigr)-f(\eta)
\bigr]^2 \bigr)+\mu \bigl(\bar{\eta}_0\bar{
\eta}_{-1} \bigl[f \bigl(\eta^{-1} \bigr)-f(\eta)
\bigr]^2 \bigr) \bigr].
\nonumber
\end{eqnarray*}
We insert the missing rates to recover terms appearing in $\mathcal
{D}_{\mathrm{FA}}(f)$. For instance, since we assumed $q<1/2$,
%
%
\begin{equation}
\mu \bigl(\bar{\eta}_0 \bigl[f \bigl(\eta^1 \bigr)-f(
\eta) \bigr]^2 \bigr)\leq\frac
{1}{q}\mu \bigl(\bar{
\eta}_0r_1(\eta) \bigl[f \bigl(\eta^1
\bigr)-f(\eta) \bigr]^2 \bigr)
\end{equation}
and
\begin{equation}
\mu \bigl(\bar{\eta}_0\bar{\eta}_1 \bigl[f
\bigl(\eta^1 \bigr)-f(\eta) \bigr]^2 \bigr)\leq
\frac{1}{p}\mu \bigl(\bar{\eta}_0r_1(\eta)
\bigl[f \bigl(\eta^1 \bigr)-f(\eta) \bigr]^2 \bigr).
\end{equation}
Thus we get
%
%
\begin{eqnarray}
\mu(jf)&\leq&\frac{q}{1-2q} \biggl\{\alpha q+\frac{1}{\alpha
q}\mathcal
{D}_{\mathrm{FA}}(f) \biggr\}+q \biggl\{\beta+\frac{1}{\beta(1-q)}\mathcal
{D}_{\mathrm{FA}}(f) \biggr\}.
\end{eqnarray}
Optimizing in $\alpha, \beta$, this yields
%
%
\begin{eqnarray}
{\mu(jf)}&\leq& \frac{2q}{1-2q}\sqrt{\mathcal{D}_{\mathrm{FA}}(f)}+2\sqrt{
\mathcal{D}_{\mathrm{FA}}(f)/(1-q)}.
\end{eqnarray}

This is enough to prove \eqref{eq:borneFA} for small $q$; see
Section~6.3 of \citet{ollalandim}.
\end{appendix}

\section*{Acknowledgments}
I am grateful to the DMA at ENS for its hospitality.
Many thanks are due to Thierry Bodineau and Cristina Toninelli
for suggesting this problem, helpful discussions and suggestions,
reading and proofreading this manuscript. I also thank an anonymous
referee for suggestions on improving the presentation.

%


\printaddresses

\end{document}